\documentclass{amsart}
\usepackage{amssymb}
\usepackage{amsmath}
\usepackage{longtable}
\newtheorem{thm}{Theorem}[section]
\newtheorem{lemma}[thm]{Lemma}
\newtheorem{definition}[thm]{Definition}

\newtheorem{corollary}[thm]{Corollary}
\newtheorem{proposition}[thm]{Proposition}

\newtheorem{notation}[thm]{Notation}

\newcommand{\ty}{\mathbf t}

\newcommand{\B}{\mathbb B}
\newcommand{\Q}{\mathbb{Q}}
\newcommand{\N}{\mathbb N}
\newcommand{\F}{\mathbb F}
\newcommand{\R}{\mathbb R}

\newcommand{\Z}{\mathbb Z}

\newcommand{\st}[1]{\vskip 1mm\noindent{\bf #1}}

\def\as#1{\renewcommand\arraystretch{#1}}
\def\bb{\mathcal{B}}

\def\dsc{\mathop{\rm Disc}}
\def\md#1{\ \mbox{\rm(mod }{#1})}
\def\gb#1{\overline{\gamma_{#1}(\t)}}

\def\j{\mathbf{j}}
\def\kb{\overline{k}}

\def\ks{k^{\op{sep}}}
\def\Ks{K^{\op{sep}}}
\def\l{\mathfrak{L}}
\def\la{\lambda}
\def\lra{\longrightarrow}
\def\lt{\ell_{\op{term}}}
\def\m{\mathfrak{m}}
\def\md#1{\ \mbox{\rm(mod }{#1})}

\def\oo{\mathcal{O}}
\def\op{\operatorname}

\def\om{\omega}
\def\ord{\operatorname{ord}}
\def\p{\mathfrak{p}}
\def\P{\mathfrak{P}}

\def\pset{\mathcal{P}}
\def\q{\mathfrak{Q}}

\def\red{\op{red}}

\def\sii{\,\Longleftrightarrow\,}
\def\st{\mathcal{S}_{\op{term}}}
\def\t{\theta}

\def\tq{\,\,|\,\,}
\def\ttt{\mathcal{T}}

\title{Higher Newton polygons and integral bases}

\thanks{Partially supported by MTM2009-13060-C02-02 and MTM2009-10359 from the Spanish MEC}

\author{Jordi Gu\`ardia}
\address{Departament de Matem\`atica Aplicada IV, Escola Polit\`ecnica Superior d'Enginyeria de Vilanova i la Geltr\'u, Av. V\'\i ctor Balaguer s/n. E-08800 Vilanova i la Geltr\'u, Catalonia, Spain}
\email{guardia@ma4.upc.edu}

\author{\hbox{Jes\'us Montes}}
\address{Departament de Ci\`encies Econ\`omiques i Empresarials,
Facultat de Ci\`encies Socials,
Universitat Abat Oliba CEU,
Bellesguard 30, E-08022 Barcelona, Catalonia, Spain\\
Departament de Matem\`atica Econ\`omica, Financera i Actuarial,
Facultat d'Economia i Empresa,
Universitat de Barcelona,
Av. Diagonal 690,
E-08034 Barcelona, Catalonia, Spain}
\email{montes3@uao.es, jesus.montes@ub.edu}

\author{\hbox{Enric Nart}}
\address{Departament de Matem\`{a}tiques,
         Universitat Aut\`{o}noma de Barcelona,
         Edifici C, E-08193 Bellaterra, Barcelona, Catalonia, Spain}
\email{nart@mat.uab.cat}

\keywords{Dedekind domain, global field, local field, Montes algorithm, Newton polygon, p-integral bases, reduced bases}

\makeatletter
\@namedef{subjclassname@2010}{%
\textup{2010} Mathematics Subject Classification}

\subjclass[2010]{Primary 11R04; Secondary 11Y40, 14G15, 14H05}

\begin{document}

\begin{abstract}
Let $A$ be a Dedekind domain whose field of fractions $K$ is a global field. Let $\p$ be a non-zero prime ideal of $A$, and $K_\p$ the completion of $K$ at $\p$. The Montes algorithm factorizes a monic irreducible separable polynomial $f(x)\in A[x]$ over $K_\p$, and it provides essential arithmetic information about the finite extensions of $K_\p$ determined by the different irreducible factors. In particular, it can be used to compute a $\p$-integral basis of the extension of $K$ determined by $f(x)$. In this paper we present a new and faster method to compute $\p$-integral bases, based    
on the use of the quotients of certain divisions with remainder of $f(x)$ that occur 
along the flow of the Montes algorithm.  
\end{abstract}
\maketitle

\section*{Introduction}Let $A$ be a Dedekind domain whose field of fractions $K$ is a global field. Let $\p$ be a non-zero prime ideal of $A$, and $\pi\in A$ a local generator of $\p$. Let  $K_\p$ be the completion of $K$ with respect to the $\p$-adic topology.

Let $f(x)\in A[x]$ be a monic irreducible separable polynomial of degree $n$. Let $\t\in\Ks$ be a root of $f(x)$, $L=K(\t)$ be the finite separable extension of $K$ generated by $\t$, and $B$ be the integral closure of $A$ in $L$.  

The Montes algorithm \cite{bordeaux,HN} computes an \emph{OM representation} of every prime ideal $\P$ of $B$ lying over $\p$ \cite{newapp}. This algorithm carries out a program suggested by \O. Ore \cite{ore23,ore28}, and developed by S. MacLane in the context of valuation theory \cite{mcla,mclb}. An OM representation is a computational object supporting several data and operators, linked to one of the irreducible factors (say) $F(x)$ of $f(x)$ in $K_\p[x]$. Among these data, the \emph{Okutsu invariants} of $F$ stand out, revealing a lot of arithmetic information about the finite extension of $K_\p$ determined by $F$ \cite{Ok,okutsu}. The initials OM stand indistinctly for Ore-MacLane or Okutsu-Montes. 

In \cite{newapp} we presented a method to compute $\p$-integral bases of $B/A$, based on these OM representations of the prime ideals of $B$ dividing $\p$. For $n$ large, this method is significantly faster than the traditional methods, most of them based on variants of the Round 2 and Round 4 routines \cite{Za,Fo,BL,VH,Ha,FPR}. 

In this paper we present an improvement of that OM-method, based on the use of \emph{quotients of $\phi$-adic expansions}. This idea goes back to a construction of integral bases by W.M. Schmidt, for certain subrings of function fields \cite{schmidt}.
Along the flow of the Montes algorithm, some polynomials $\phi(x)\in A[x]$ are constructed as a kind of optimal approximations to the irreducible factors of $f(x)$ over $K_\p$. The (conveniently truncated) $\phi$-expansions of $f(x)$ provide the necessary data to build  higher order Newton polygons of $f(x)$, from which new and better approximations are deduced. 
As a by-product of the computation of any $\phi$-expansion, $f(x)=\sum_{0\le s}a_s(x)\phi(x)^s$, we obtain several quotients: 
$$
f(x)=\phi(x)Q_1(x)+a_0(x),\quad
Q_1(x)=\phi(x)Q_2(x)+a_1(x),\quad\dots
$$
These polynomials $Q_i(x)$ have nice properties that can be exploited to obtain shortcuts and improvements in the computation of $\p$-integral bases. 

The outline of the paper is as follows. In section \ref{secOkutsu} we review the main technical ingredients of the paper: OM representations and Okutsu invariants of irreducible separable polynomials over local fields. In section \ref{secOM} we review the OM-method of \cite{newapp} for the computation of integral bases. In section \ref{secMain}, we study the quotients $Q(x)$ obtained along the computation of $\phi$-expansions
of $f(x)$. We analyze the $\P$-adic value of $Q(\t)$, for all prime ideals $\P$ lying over $\p$, in order to determine the highest exponent $\mu$ such that $Q(\t)/\pi^{\mu}$ is $\p$-integral (Theorem \ref{denquot} and Corollary \ref{applications}). In section \ref{secIB}, we show how to construct local bases with these quotients. 
For every prime ideal $\P\mid \p$, we find a family of elements of $L$ whose images in the $\P$-completion $L_\P$ are an integral basis of the local extension $L_\P/K_\p$. These elements are constructed
as a product of quotients, divided by an adequate power of $\pi$. The essential difference with the OM-method is that all these elements are already $\p$-integral (globally integral if $A$ is a PID), and not only $\P$-integral. Finally, in section \ref{secMethod} we show how to use these $\p$-integral elements to build a $\p$-integral basis 
(Theorem \ref{pBasis2}). This \emph{method of the quotients}  has three significant advantages with respect to the former OM-method and all classical methods:\medskip

(1) \ It yields always \emph{$\p$-reduced} bases. For instance, let $L=\F(t,x)$ be the function field of a curve $C$ over a finite field $\F$, defined by an equation $f(t,x)=0$, which is separable over $K=\F(t)$. For the subring $A=\F[t^{-1}]$ and the prime ideal $\p=t^{-1}A$, a $\p$-reduced basis of $B/A$ is just a classical reduced basis with respect to a certain size function determined by the degree function on $A$ \cite{aklenstra}, \cite[Sec. 16]{hwlenstra}. The construction of reduced bases is a key ingredient in the computation of bases of the Riemann-Roch spaces attached to divisors of $C$ \cite{schmidt}, \cite{schoernig}, \cite{hess}. \medskip

(2) \ It admits a neat complexity analysis (Theorem \ref{complexity}).   
The method requires only $O(n)$ multiplications in the ring $A[\t]$, along an ordinary application of the Montes algorithm with input data $(f(x),\p)$. If $A/\p$ is  small, the computation of a $\p$-integral basis requires altogether $O\left(n^{2+\epsilon}\delta^{1+\epsilon}+n^{1+\epsilon}\delta^{2+\epsilon}\right)$ word operations, where $\delta$ is the $\p$-adic valuation of the discriminant of $f(x)$.  \medskip

(3) \ It has an excellent practical performance. For $A=\Z$, the method may be tested by running the {\tt pIntegralBasis} routine of the Magma package {\tt +Ideals.m}, which may be downloaded from the site
{\tt http://www-ma4.upc.edu/$\sim$guardia/} {\tt +Ideals.html}.\medskip


\section{Okutsu invariants of irreducible polynomials over local fields}\label{secOkutsu} 
Let $k$ be a local field, i.e. a locally compact and complete field with respect to a discrete valuation $v$. Let $\oo$ be the valuation ring of $k$, $\m$ the maximal ideal, $\pi\in\m$ a generator of $\m$ and $\F=\oo/\m$ the residue field, which is a finite field. 

Let $\ks\subset \kb$ be the separable closure of $k$ inside a fixed algebraic closure. Let $v\colon \kb\to \Q\cup\{\infty\}$, be the canonical extension of the discrete valuation $v$ to $\kb$, normalized by $v(k)=\Z$.

Let $F(x)\in\oo[x]$ be a monic irreducible separable polynomial, $\t\in \ks$ a root of $F(x)$, and $L=k(\t)$ the finite separable extension of $k$ generated by $\t$. Denote $n:=[L\colon k]=\deg F$. Let $\oo_L$ be the ring of integers of $L$, $\m_L$ the maximal ideal and $\F_L$ the residue field. 
We indicate with a bar, $\raise.8ex\hbox{---}\colon \oo[x]\longrightarrow \F[x]$,
the canonical homomorphism of reduction of polynomials modulo $\m$. 

Let $[\phi_1,\dots,\phi_r]$ be an \emph{Okutsu frame} of $F(x)$, and let $\phi_{r+1}$ be an \emph{Okutsu approxi\-mation} to $F(x)$. That is, $\phi_1,\dots,\phi_{r+1}\in\oo[x]$ are monic separable polynomials of strictly increasing degree: $$1\le m_1:=\deg\phi_1<\cdots <m_r:=\deg\phi_r<m_{r+1}:=\deg\phi_{r+1}=n,$$ and for any monic polynomial $g(x)\in\oo[x]$ we have: 
\begin{equation}\label{frame}
m_i\le \deg g<m_{i+1}\ \Longrightarrow\ \dfrac{v(g(\t))}{\deg g}\le\dfrac{v(\phi_i(\t))}{m_i}<\dfrac{v(\phi_{i+1}(\t))}{m_{i+1}},
\end{equation}
for $0\le i\le r$, with the convention that $m_0=1$ and $\phi_0(x)=1$. It is easy to deduce from (\ref{frame}) that the polynomials $\phi_1(x),\dots,\phi_{r+1}(x)$ are all irreducible in $\oo[x]$. 

The length $r$ of the frame is called the \emph{Okutsu depth} of $F(x)$. We have $r=0$ if and only if $\overline{F}$ is irreducible over $\F$; in this case, the Okutsu frame is an empty list.
Okutsu frames were introduced by K. Okutsu in \cite{Ok} as a tool to construct integral bases. Okutsu approximations were introduced in \cite{okutsu}, where it is shown that the family $\phi_1,\dots,\phi_{r+1}$ determines an \emph{optimal $F$-complete type of order $r+1$}:
\begin{equation}\label{OM}
\ty_F=
(\psi_0;(\phi_1,\lambda_1,\psi_1);\cdots;(\phi_r,\lambda_r,\psi_r);(\phi_{r+1},\lambda_{r+1},\psi_{r+1})).
\end{equation}
In the special case $\phi_{r+1}=F$, we have $\la_{r+1}=-\infty$ and $\psi_{r+1}$ is not defined. We call $\ty_F$ an \emph{OM representation} of $F$.

Any OM representation of the polynomial $F$ carries (stores) several invariants and operators yielding strong arithmetic information about $F$ and the extension $L/k$.
Let us recall some of these invariants and operators.

Attached to the type $\ty_F$, there is a family of discrete valuations of the rational function field $k(x)$, the \emph{MacLane valuations}:
$$
v_i\colon k(x)\longrightarrow \Z\cup\{\infty\},\quad 1\le i\le r+1,
$$
satisfying $0=v_1(F)<\cdots<v_{r+1}(F)$. The $v_1$-value of a polynomial in $k[x]$ is the minimum of the $v$-values of its coefficients.

Also, $\ty_F$ determines a family of Newton polygon operators:
$$
N_i\colon k[x]\longrightarrow 2^{\R^2},\quad 1\le i\le r+1,
$$where $2^{\R^2}$ is the set of subsets of the Euclidean plane. Any non-zero polynomial $g(x)\in k[x]$ has a canonical $\phi_i$-development:
$$
g(x)=\sum\nolimits_{0\le s}a_s(x)\phi_i(x)^s,\quad \deg a_s<m_i,
$$ and the polygon $N_i(g)$ is the lower convex hull of the set of points $(s,v_i(a_s\phi_i^s))$. Usually, we are only interested in the principal polygon $N_i^-(g)\subset N_i(g)$ formed by the sides of negative slope. For all $1\le i\le r$, the Newton polygons $N_i(F)$ and $N_i(\phi_{i+1})$ are one-sided and they have the same slope, which is a negative rational number $\lambda_i\in\Q_{<0}$.  
The Newton polygon $N_{r+1}(F)$ is one-sided and it has an (extended) integer negative slope, which we denote by $\lambda_{r+1}\in\Z_{<0}\cup\{-\infty\}$.  

There is a chain of finite extensions: $\F=\F_0\subset \F_1\subset\cdots\subset\F_{r+1}=\F_L$. The type $\ty_F$ stores monic irreducible polynomials $\psi_i(y)\in\F_i[y]$ such that $\F_{i+1}\simeq \F_i[y]/(\psi_i(y))$. We have $\psi_i(y)\ne y$, for all $i>0$. Finally, for every negative rational number $\lambda$, there are \emph{residual polynomial} operators:
$$
R_{\lambda,i}\colon k[x] \longrightarrow \F_i[y],\quad 0\le i\le r+1.
$$ 
We define $R_i:=R_{\lambda_i,i}$. For all $0\le i\le r$, we have $R_i(F)\sim \psi_i^{\omega_{i+1}}$ and $R_i(\phi_{i+1})\sim\psi_i$, where the symbol $\sim$ indicates that the polynomials coincide up to a multiplicative constant in $\F_i^*$. For $i=0$ we have $R_0(F)=\overline{F}=\psi_0^{\omega_1}$ and $R_0(\phi_1)=\overline{\phi_1}=\psi_0$. The exponents $\omega_{i+1}$ are all positive and $\omega_{r+1}=1$. The operator $R_{r+1}$ is defined only when $\phi_{r+1}\ne F$; in this case, we also have $R_{r+1}(F)\sim \psi_{r+1}$, with $\psi_{r+1}(y)\in\F_{r+1}[y]$ monic of degree one such that $\psi_{r+1}(y)\ne y$. 

From these data some more numerical invariants are deduced. Initially we take:
$$m_0:=1,\quad f_0:=\deg \psi_0,\quad e_0:=1,\quad h_0:=V_0=0.
$$Then, we define for all $1\le i\le r+1$:
$$
\as{1.2}
\begin{array}{l}
h_i,\,e_i \ \mbox{ positive coprime integers such that }\la_i=-h_i/e_i,\\
f_i:=\deg \psi_i,\\
m_i:=\deg \phi_i=e_{i-1}f_{i-1}m_{i-1}=(e_0\,e_1\cdots e_{i-1})(f_0f_1\cdots f_{i-1}),\\
V_i:=v_i(\phi_i)=e_{i-1}f_{i-1}(e_{i-1}V_{i-1}+h_{i-1}),\\
\ell_i,\,\ell'_i \ \mbox{ a pair of integers such that }\ell_ih_i-\ell'_ie_i=1,\\
z_{i-1}:= \ \mbox{the class of $y$ in $\F_{i}$, so that }\psi_{i-1}(z_{i-1})=0.
\end{array}
$$
\as{1.}

An irreducible polynomial $F$ admits infinitely many different OM representations. However, the numeri\-cal invariants $e_i,f_i,h_i$, for $0\le i\le r$, and the MacLane valuations 
$v_1,\dots,v_{r+1}$ attached to $\ty_F$, are canonical invariants of $F$. 

The data $\lambda_{r+1},\psi_{r+1}$ are not invariants of $F$; they depend on the choice of the Okutsu approximation $\phi_{r+1}$. The integer slope $\lambda_{r+1}=-h_{r+1}$ measures how close is $\phi_{r+1}$ to $F$. We have $\phi_{r+1}=F$ if and only if $h_{r+1}=\infty$.
 
\begin{definition}
An \emph{Okutsu invariant} of $F(x)$ is a rational number that depends only on $e_0,e_1,\dots,e_r,f_0,f_1,\dots,f_r,h_1,\dots,h_r$.
\end{definition}

For instance, the ramification index and residual degree of $L/k$ are Okutsu invariants of $F$. More precisely,
$$
e(L/k)=e_0e_1\cdots e_r,\quad f(L/k)=f_0f_1\cdots f_r.
$$

The general definition of a \emph{type} may be found in \cite[Sec. 2.1]{HN}. In later sections, we shall consider types which are not necessarily optimal nor $F$-complete. So, it may be convenient to distinguish these two properties among all features of a type that we have just described.

\begin{definition}\label{optimal}Let  $\ty=(\psi_0;(\phi_1,\lambda_1,\psi_1);\cdots;(\phi_i,\lambda_i,\psi_i))$ be a type of order $i$ and denote $m_{i+1}:=e_if_im_i$. Let $g(x),h(x)\in k[x]$ be non-zero polynomials.\medskip

\noindent{$\bullet$} We say that $\ty$ is \emph{optimal} if $m_1<\cdots<m_i$.
We say that $\ty$ is \emph{strongly optimal} if $m_1<\cdots<m_i<m_{i+1}$.
\medskip

\noindent{$\bullet$} We define $\ord_\ty(g):=\ord_{\psi_{i}}R_{i}(g)$ in $\F_i[y]$.
If $\ord_\ty(g)>0$, we say that $\ty$ \emph{divides} $g(x)$, and we write $\ty\mid g(x)$. We have $\ord_\ty(gh)=\ord_\ty(g)+\ord_\ty(h)$.  
\medskip

\noindent{$\bullet$}We say that $\ty$ is \emph{$g$-complete} if $\ord_\ty(g)=1$.
\medskip

\noindent{$\bullet$} A \emph{representative} of $\ty$ is a monic polynomial $\phi(x)\in\oo[x]$ of degree $m_{i+1}$, such that $R_i(\phi)\sim \psi_i$. This polynomial is necessarily irreducible in $\oo[x]$. The degree $m_{i+1}$ is minimal among all polynomials satisfying this condition. 
The choice of a representative of $\ty$ determines a Newton polygon operator $N_{\phi,v_{i+1}}$, which depends only on $\phi$ and the valuation $v_{i+1}$ supported by $\ty$. We denote $N_{i+1}:=N_{\phi_{i+1},v_{i+1}}$.\medskip

\noindent{$\bullet$} For any $0\le j\le i$, the \emph{truncation} of $\ty$ at level $j$, $\op{Trunc}_j(\ty)$, is the type of order $j$ obtained from $\ty$ by dropping all levels  higher than $j$.  
\end{definition}

For a general type of order $i$ dividing $F$, we have $m_1\mid\cdots\mid m_i$ and $\omega_i>0$, but not necessarily $m_1<\cdots<m_i=\deg F$, and $\omega_i=1$. These were particular properties of the optimal and $F$-complete type $\ty_F$ of order $i=r+1$, constructed from an Okutsu frame and an Okutsu approximation to $F$. \medskip

\begin{definition}
The \emph{length} of a Newton polygon $N$ is the abscissa of its right end point; we denote it by $\ell(N)$. 
\end{definition}

\begin{lemma}{\cite[Lem. 2.17,(2)]{HN}}\label{length}
Let $\ty$ be a type of order $i\ge 0$, and let $\phi_{i+1}\in\oo[x]$ be a representative of $\ty$. 
Then, $\ell(N_{i+1}^-(g))=\ord_\ty(g)$, for any non-zero polynomial $g(x)\in k[x]$.
\end{lemma}

The next lemma shows that a type gathers ``features" of irreducible polynomials in $\oo[x]$.  This is the motivation for the term \emph{type}. The lemma is a combination of \cite[Def. 2.1, Lem. 2.4, Cor. 2.18]{HN}.

\begin{lemma}\label{type}
Let $\ty$ be a type of order $i\ge 0$, and let $G(x)\in\oo[x]$ be a monic irreducible separable polynomial such that $\ty\mid G$. Then,
\begin{enumerate}
\item $R_j(G)\sim\psi_j^{n_j}$ in $\F_j[y]$, for a certain $n_j>0$, for all $0\le j\le i$.
\item $N_j(G)=N_j^-(G)$ is one-sided of slope $\la_j$, for all $1\le j\le i$.
\end{enumerate}
\end{lemma}

We shall frequently use the following result, extracted from \cite[Prop. 3.5,(5)]{HN}.
Note that it contains the definition of the MacLane valuation $v_{i+1}$. 

\begin{proposition}\label{vgt}
Let $\ty$ be a type of order $i\ge 1$, and let $F(x)\in\oo[x]$ be a monic irreducible separable polynomial such that $\ty\mid F$. Let $\t\in\ks$ be a root of $F(x)$, and
 $g(x)\in\oo[x]$ a non-zero polynomial. Take a line of slope $ \lambda_i$ far below $N_i(g)$, and shift it upwards till it touches the polygon for the first time. Let $(0,H)$ be the intersection point of this line with the vertical axis. Then,  
$$v(g(\t))\ge v_{i+1}(g)/(e_0\cdots e_i)=H/(e_0\cdots e_{i-1}),$$ and equality holds if and only if $\ty\nmid g(x)$.
\end{proposition}

\begin{corollary}\label{previous}
With the above notation,  
$v(\phi_j(\t))=(V_j+|\lambda_j|)/(e_0\cdots e_{j-1})$, for all $1\le j\le i$.
\end{corollary}

\subsection*{Local integral bases}
The next result is an elementary combinatorial fact.

\begin{lemma}\label{combasic}For $0\le i\le r$,  consider positive integers, $m_i\mid\cdots\mid m_{r+1}$. Then, any integer $0\le N <(m_{r+1}/m_i)$ can be expressed in a unique way as:
$$
N=j_i+j_{i+1}(m_{i+1}/m_i)+\cdots+j_r(m_r/m_i),
$$
for integers $j_k$ satisfying: $0\le j_k<(m_{k+1}/m_k)$, for all $i\le k\le r$.
\end{lemma}

Now, let $m_1,\dots,m_{r+1}$ be the degrees of the Okutsu polynomials of an OM representation $\ty_F$ of a monic irreducible separable polynomial $F\in\oo[x]$, as in (\ref{OM}). Take $m_0:=1$, $\phi_0(x):=x$. For any integer $0\le m<m_{r+1}=n$, consider the following monic polynomial $g_m(x)\in\oo[x]$, of degree $m$:
$$
g_m(x):=\prod\nolimits_{0\le k\le r}\phi_k(x)^{j_k},\quad m=\sum\nolimits_{0\le k\le r}j_km_k, \ 0\le j_k<m_{k+1}/m_k=e_kf_k.
$$
Corollary \ref{previous} provides concrete formulas for all $v(\phi_k(\t))$; thus, we can easily compute $\mu_m:=\lfloor v(g_m(\t))\rfloor $, for all $m$.

\begin{thm}[Okutsu, {\cite[I,Thm. 1]{Ok}}]\label{basis}
The following family is an $\oo$-basis of $\oo_L$:
$$
1,\ g_1(\t)/\pi^{\mu_1},\ \dots,\ g_{n-1}(\t)/\pi^{\mu_{n-1}}.
$$
\end{thm}

The \emph{exponent} of $F$ is the least non-negative integer $\exp(F)$ such that $\pi^{\exp(F)}\oo_L$ is included in $\oo[\t]$.
Since $\mu_1\le \cdots\le\mu_{n-1}$, it is clear that
$\exp(F)=\mu_{n-1}$.

\section{An OM method to compute $\mathbf{\p}$-integral bases}\label{secOM}

Let $A$ be a Dedekind domain whose field of fractions $K$ is a global field, and let $\Ks$ be a separable closure of $K$. Let $f(x)\in A[x]$ be a monic irreducible and separable polynomial of degree $n>1$. Let $L=K(\t)$ be the finite separable extension of $K$ generated by a root $\t\in\Ks$ of $f(x)$. The integral closure $B\subset L$ of $A$ in $L$ is a Dedekind domain too.

Let $\p$ be a non-zero prime ideal of $A$. Let 
$A_\p$ be the localization of $A$ at $\p$, $\pi\in A$ a generator of the principal ideal $\p A_\p$, and $\F_\p:=A/\p$ the residue field.
 The integral closure $B_\p$ of $A_\p$ in $L$ is the subring of $\p$-integral elements of $L$:
$$
B_\p=\{\alpha\in L\mid v_\P(\alpha)\ge 0,\ \forall \P\in \op{Spec}(B),\ \P\mid \p\},
$$ where $v_\P$ is the discrete valuation of $L$ attached to $\P$.
The ring $B_\p$ is a free $A_\p$-module of rank $n$. 

\begin{definition}\label{pintegral}
A $\p$-integral basis of $B/A$ is a family $\alpha_1,\dots,\alpha_n\in B_\p$, that satisfies any of the following equivalent conditions:

\begin{enumerate}
\item[(a)]  $\alpha_1,\dots,\alpha_n$ is an $A_\p$-basis of $B_\p$.
\item[(b)] $\alpha_1\otimes 1,\dots,\alpha_n\otimes 1$ is an $\F_\p$-basis of $B_\p\otimes_{A_\p}\F_\p\simeq B_\p/\p B_\p\simeq B/\p B$.
\end{enumerate}
\end{definition}

Conditions (a) and (b) are equivalent by Nakayama's lemma.
Since $B/\p B$ has dimension $n$ as an $\F_\p$-vector space, it suffices to check that $\alpha_1,\dots,\alpha_n\in B_\p$ determine $\F_\p$-linearly independent elements in the $\F_\p$-algebra $B/\p B$, to show that they form a $\p$-integral basis of $B/A$. 

Consider the factorization of $\p B$ into a product of prime ideals in $L$:
$$
\p B=\P_1^{e(\P_1/\p)}\cdots\, \P_g^{e(\P_g/\p)}.
$$
Let $K_\p$, $L_\P$, be the completions of $K$ and $L$ with respect to the $\p$-adic and $\P$-adic topology, respectively. Denote the ring of integers of these fields by:
$$
\mathcal{O}_\p\subset K_\p,\qquad \oo_\P\subset L_\P,\ \forall\,\P\mid \p.
$$
Finally, we denote by $n_\P:=[L_\P\colon K_\p]=e(\P/\p)f(\P/\p)$, the local degrees.

By a classical theorem of Hensel, these prime ideals are in 1-1 correspondence with the different monic irreducible factors of $f(x)$ in $\oo_\p[x]$.

\begin{definition}
For each prime ideal $\P\mid \p$, let us fix a topological embedding, $i_\P\colon L\subset L_\P\hookrightarrow \overline{K}_\p$. Then $\t_\P:=i_\P(\t)$ is the root of a unique monic irreducible factor (say) $F_\P(x)$ of $f(x)$ over $\oo_\p$. Also, we denote:
$$w_\P:=e(\P/\p)^{-1}v_\P\colon L^*\lra e(\P/\p)^{-1}\Z.$$
\end{definition}

Clearly, $w_\P(\alpha)=v(i_\P(\alpha))$, for all $\alpha\in L$, where $v:=v_\p$ is the canonical extension of $v_\p$ to $\overline{K}_\p$.
Thus, for any polynomial $g(x)\in A[x]$,
$$
w_\P(g(\t))=v(g(\t_\P)).
$$
This identity will be implicitly used throughout the paper without further mention, when we apply local results to a global situation.

The Montes algorithm provides a family of OM representations of the irreducible factors of $f(x)$ in $\oo_\p[x]$. For any prime ideal $\P$ dividing $\p$, let us denote by 
$$
\ty_\P:=\ty_{F_\P}=\left(\psi_{0,\P};(\phi_{1,\P},\lambda_{1,\P},\psi_{1,\P});\cdots;(\phi_{r_\P+1,\P},\lambda_{r_\P+1,\P},\psi_{r_\P+1,\P})\right),
$$
the OM representation corresponding to $F_\P$. All polynomials $\phi_{i,\P}$ have coefficients in $A$. The type $\ty_\P$ singles out $\P$ (or $F_\P$) by:
$$
\ty_\P\mid F_\P,\quad \ty_\P\nmid F_\q, \quad \forall\,\q\mid\P,\ \q\ne\P.
$$

As we saw in the last section, from the OM representation we derive a family of $\P$-integral elements in $L$,
$$
\bb_\P=\left\{ 1,\ g_{1,\P}(\t)/\pi^{\mu_{1,\P}},\ \dots,\ g_{n_\P-1,\P}(\t)/\pi^{\mu_{n_\P-1,\P}}\right\}\subset L.
$$
whose image under $i_\P$ is the Okutsu basis of $\oo_\P$ as an $\oo_\p$-module, described in Theorem \ref{basis}.
It is easy to buid a $\p$-integral basis with these local $\P$-bases.

\begin{thm}[Ore, \cite{ore25}]\label{ore}
For each prime ideal $\P\mid\p$, take $\beta_\P\in B_\p$ such that:
$$
w_\P(\beta_\P)=0,\quad w_\q(\beta_\P)\ge \exp(F_\P)+1, \ \forall\,\q\mid \p,\ \q\ne\P.
$$
Then, $\bb:=\bigcup_{\P\mid \p}\beta_\P\bb_\P$, 
is a $\p$-integral basis.
\end{thm}

In \cite[Secs. 3.2,4.2]{newapp} we found an efficient way to compute these multiplicators $\beta_\P$ in terms of the data supported by the OM representations. They are elements in $B_\p$ of the form:
$$
\beta_\P=\pi^{-N}\prod\nolimits_{\q\mid\p,\,\q\ne\P}\phi_\q(\t)^{d_\q}, \quad \phi_\q:=\phi_{r_\q+1,\q}
$$ 
An adequate choice of the exponents $d_\q,\,N$ leads to $w_\P(\beta_\P)=0$. Also, we can get $w_\q(\beta_\P)$ high enough, by improving to an adequate precision the Okutsu approximations $\phi_\q$ to the factors $F_\q$, with the single-factor lifting algorithm \cite{GNP}.

Although the local bases are a by-product of the Montes algorithm, the cons\-truction of these multiplicators requires an extra work. In section \ref{secIB} we shall cons\-truct local bases by using \emph{quotients} instead of $\phi$-polynomials. These quotients are $\p$-integral (and not only $\P$-integral). In section \ref{secMethod} we shall use this fact to construct $\p$-integral bases with no need to compute multiplicators (Theorem \ref{pBasis2}).   

If $A$ is a PID, then $B$ is a free $A$-module of rank $n$. If we take $\pi\in A$ to be a generator of the principal ideal $\p$, then the $\p$-integral bases constructed as above are made of global integral elements, because $\p$ is the only prime ideal of $A$ that divides the denominators.

If for all prime ideals $\p$ dividing the discriminant of $f(x)$ we compute a $\p$-integral basis in Hermite normal form, then an easy application of the CRT yields a global integral basis; that is, a basis of $B$ as an $A$-module.

\section{Quotients of $\phi$-adic expansions}\label{secMain}
We keep all notation from the preceding section.

Let $\phi(x)\in A[x]$ be a monic polynomial of positive degree, and let
$$
f(x)=a_0(x)+a_1(x)\phi(x)+\cdots+a_m(x)\phi(x)^m,\quad a_s(x)\in A[x],\ \deg a_s<\deg \phi,
$$ be the canonical $\phi$-expansion of $f(x)$. Note that $m=\lfloor \deg f/\deg \phi\rfloor$.

\begin{definition}
The \emph{$\phi$-quotients of $f(x)$} are the quotients  $Q_1(x),\dots,Q_m(x)$, obtained along the computation of the coefficients of the $\phi$-expansion of $f(x)$:
$$
\begin{array}{rl}
f(x)=\!\!&\!\phi(x)\,Q_1(x)+a_0(x),\\
Q_1(x)=\!\!&\!\phi(x)\,Q_2(x)+a_1(x),\\
\cdots&\cdots\\
Q_m(x)=\!\!&\!\phi(x)\cdot 0+a_m(x)=a_m(x).
\end{array}
$$
\end{definition}

Equivalently, $Q_s(x)$ is the quotient of the division of $f(x)$ by $\phi(x)^s$; we denote by $r_s(x)$ the remainder of this division. Thus, for all $1\le i\le m$ we have,
\begin{equation}\label{residue}
 f(x)=r_s(x)+Q_s(x)\phi(x)^s,\ r_s(x)=a_0(x)+a_1(x)\phi(x)+\cdots+a_{s-1}(x)\phi(x)^{s-1}.
\end{equation}

The aim of this section is to use these quotients to construct nice $\p$-integral elements. More precisely, for certain $\phi$-quotients $Q(x)$ of $f(x)$, we find the highest exponent $\mu$ such that $Q(\t)/\pi^\mu$ is $\p$-integral (Theorem \ref{denquot} and Corollary \ref{applications}). 

\subsection{Construction of integral elements}\label{subsecDenQ}
\begin{lemma}\label{multiadic}
Let $\ty=(\psi_0;\cdots;(\phi_{r},\lambda_{r},\psi_{r}))$ be a type over $\oo_\p$, of order $r\ge1$. Fix an index $1\le i\le r$. Let $g(x)\in A[x]$ be a polynomial of degree less than $m_{r+1}$, and consider its multiadic expansion:
$$
g(x)=\sum\nolimits_{\j=(j_i,\dots,j_r)}a_{\j}(x)\Phi(x)^{\j}, \quad\deg a_{\j}<m_i,
$$
where $\Phi(x)^{\j}:=\phi_i(x)^{j_i}\cdots \phi_r(x)^{j_r}$, and $0\le j_k<e_kf_k$, for all $i\le k\le r$. Then, $v_{r+1}(g)=\min\left\{v_{r+1}(a_{\j}(x)\Phi(x)^\j)\mid \j=(j_i,\dots,j_r)\right\}$.
\end{lemma}

\begin{proof}
Since $v_{r+1}$ is a valuation, it is sufficient to show $v_{r+1}\left(a_{\j}(x)\Phi(x)^\j\right)\ge v_{r+1}(g)$, for all $\j$. Let us prove this inequality by induction on $r-i$. For $r=i$ this is proven in \cite[Prop. 2.7,(4)]{HN}. Suppose that
$r>i$ and the lemma is true for the indices $r-1\ge i$. 
Consider the $\phi_r$-expansion of $g(x)$, and the $(\phi_i,\dots,\phi_{r-1})$-multiadic development of each coefficient:
$$
g(x)=\sum_{0\le j<e_rf_r}g_j(x)\phi_r(x)^j,\quad 
g_j(x)=\sum_{\j=(j_i,\dots,j_{r-1},j)}a_\j(x)\phi_i(x)^{j_i}\cdots \phi_{r-1}(x)^{j_{r-1}}.
$$
By the definition of $v_{r+1}$, we have 
$v_{r+1}(P)=e_rv_r(P)$, for any polynomial $P\in A[x]$ of degree less than $m_r$. Thus, by \cite[Prop. 2.7,(4)]{HN} and the induction hypothesis:
\begin{align*}
v_{r+1}(g)&\le\, v_{r+1}\left(g_j(x)\phi_r(x)^j\right)=e_rv_r(g_j)+jv_{r+1}(\phi_r)\\
&\le\, e_r v_r\left(a_\j(x)\phi_i(x)^{j_i}\cdots \phi_{r-1}(x)^{j_{r-1}}\right)+jv_{r+1}(\phi_r)=\,v_{r+1}\left(a_\j(x)\Phi(x)^\j\right),
\end{align*}
for all $0\le j<e_rf_r$, and all $\j=(j_i,\dots,j_r)$ such that $j_r=j$. 
\end{proof}

For any pair $i<r$ of positive integers, two of the formulas from \cite[Prop. 2.15]{HN} can be rewritten as:
\begin{equation}\label{vrphiHN}
\dfrac{v_r(\phi_r)}{e_1\cdots e_{r-1}}=\sum_{1\le j<r}\dfrac{m_r}{m_j}\,\dfrac{h_j}{e_1\cdots e_j},\quad 
\dfrac{v_r(\phi_i)}{e_1\cdots e_{r-1}}=\sum_{1\le j\le i}\dfrac{m_i}{m_j}\,\dfrac{h_j}{e_1\cdots e_j}.
\end{equation}
Recall that $V_i:=v_i(\phi_i)$. We deduce from these identities:
\begin{equation}\label{vrphir}
\dfrac{V_r}{e_1\cdots e_{r-1}}=\dfrac{m_r}{m_i}\,\dfrac{V_i}{e_0\cdots e_{i-1}}+\sum_{i\le j<r}\dfrac{m_r}{m_j}\,\dfrac{h_j}{e_1\cdots e_j},
\end{equation}
\begin{equation}\label{vrphis}
\dfrac{v_r(\phi_i)}{e_1\cdots e_{r-1}}=\dfrac{V_i}{e_0\cdots e_{i-1}}+\dfrac{h_i}{e_1\cdots e_{i}}.
\end{equation}

\begin{thm}\label{denquot}
Let $\ty=(\psi_0;\cdots;(\phi_{r-1},\lambda_{r-1},\psi_{r-1}))$ be a type over $\oo_\p$, of order $r-1\ge0$, and let $\phi_r$ be a representative of $\ty$. Suppose that $\ty\mid f(x)$ and all po\-lynomials $\phi_1,\dots,\phi_{r}$ have coefficients in $A$. For any integer, $1\le s\le \ell(N_r^-(f))$, let $Q_s(x)$ be the $s$-th $\phi_{r}$-quotient of $f(x)$, and let $y_s\in\Q$ be determined by $(s,y_s)\in N_{r}^-(f)$. Then, for every prime ideal $\P$ of $B$ lying over $\p$, we have:
\begin{equation}\label{aim}
w_{\P}(Q_s(\t))\ge H_s:=(y_s-sV_{r})/(e_0\cdots e_{r-1}).
\end{equation}
In particular, $Q_s(\t)/\pi^{\lfloor H_s\rfloor}$ is $\p$-integral. 
\end{thm}

\begin{proof}
Let $\la_{r}=-h_{r}/e_{r}$, with $h_{r},\,e_{r}$ positive coprime integers, be the slope of the side $S$ of $N_{r}^-(f)$, whose projection to the horizontal axis contains the abscissa $s$. If $s$ is the abscissa of a vertex of $N_{r}^-(f)$, then we take $S$ to be the left adjacent side. 

The identities (\ref{residue}) show that $N_{r}^-(f)$ should split in principle into two parts: $N_{r}^-(r_s)$ and $N_{r}^-(Q_s(\phi_{r})^s)$ (see Figure \ref{figSplit}).
This is not always true because, depending on the values of $v_{r}(a_{s-1}(\phi_{r})^{s-1})$ and $v_{r}(a_s(\phi_{r})^s)$, the two parts of the side $S$ in the polygons $N_{r}^-(r_s)$ and
$N_{r}^-(Q_s(\phi_{r})^s)$  might change. Figure \ref{figSplit2} shows different possibilities for these changes. However, the line $L$ of slope $\la_{r}$ that first touches both polygons from below is still the line determined by $S$.

\begin{figure}\caption{}\label{figSplit}
\setlength{\unitlength}{5.mm}
\begin{picture}(11,10)
\put(-.15,9.85){$\bullet$}\put(1.85,5.85){$\bullet$}
\put(5.85,3.85){$\bullet$}\put(10.85,2.85){$\bullet$}
\put(0,-.5){\line(0,1){11}}\put(-1,1){\line(1,0){13}}
\put(2,6.03){\line(-1,2){2}}\put(2,6){\line(-1,2){2}}
\put(6,4){\line(-2,1){4}}\put(6,4.03){\line(-2,1){4}}
\put(11,3){\line(-5,1){5}}\put(11,3.03){\line(-5,1){5}}
\put(3.75,4.85){$\times$}\put(4.75,4.35){$\times$}
\put(5.5,7){\begin{footnotesize}$N_{r}^-(f)$\end{footnotesize}}
\put(3.2,5.6){\begin{footnotesize}$S$\end{footnotesize}}
\multiput(4,.9)(0,.25){17}{\vrule height2pt}
\multiput(5,.9)(0,.25){15}{\vrule height2pt}
\multiput(11,.9)(0,.25){9}{\vrule height2pt}
\put(2.8,1.2){\begin{footnotesize}$s\!-\!1$\end{footnotesize}}
\put(5.2,1.2){\begin{footnotesize}$s$\end{footnotesize}}
\multiput(-.1,3)(.25,0){45}{\hbox to 2pt{\hrulefill }}
\put(-1.6,2.9){\begin{footnotesize}$v_{r}(f)$\end{footnotesize}}
\put(.2,1.2){\begin{footnotesize}$0$\end{footnotesize}}
\put(2,.5){\vector(-1,0){2}}\put(2,.5){\vector(1,0){2}}
\put(7,.5){\vector(-1,0){2}}\put(7,.5){\vector(1,0){4}}
\put(1,-.2){\begin{footnotesize}$N_{r}^-(r_s)$\end{footnotesize}}
\put(6.2,-.2){\begin{footnotesize}$N_{r}^-(Q_s(\phi_{r})^s)$\end{footnotesize}}
\put(-.8,4.4){\begin{footnotesize}$y_s$\end{footnotesize}}
\multiput(-.1,4.5)(.25,0){21}{\hbox to 2pt{\hrulefill }}
\end{picture}
\end{figure}
\begin{figure}\caption{}\label{figSplit2}
\setlength{\unitlength}{5.mm}
\begin{picture}(11,10)
\put(-.15,9.85){$\bullet$}\put(1.85,5.85){$\bullet$}\put(5.85,3.85){$\bullet$}\put(10.85,2.85){$\bullet$}
\put(3.85,7.85){$\bullet$}\put(2.85,5.85){$\bullet$}\put(4.85,6.85){$\bullet$}
\put(0,-.5){\line(0,1){11}}\put(-1,1){\line(1,0){13}}
\put(2,6.03){\line(-1,2){2}}\put(2,6){\line(-1,2){2}}
\put(11,3){\line(-5,1){5}}\put(11,3.03){\line(-5,1){5}}
\put(2,6){\line(1,0){1}}\put(2,6.03){\line(1,0){1}}
\put(3,6){\line(1,2){1}}\put(3,6.03){\line(1,2){1}}
\put(6,4){\line(-1,3){1}}\put(6,4.03){\line(-1,3){1}}
\put(2,6){\line(2,-1){2}}\put(2,6.03){\line(2,-1){2}}
\put(5,4.5){\line(2,-1){1}}\put(5,4.53){\line(2,-1){1}}
\put(3.85,4.85){$\bullet$}\put(4.85,4.35){$\bullet$}
\multiput(4,.9)(0,.25){28}{\vrule height2pt}
\multiput(5,.9)(0,.25){25}{\vrule height2pt}
\multiput(11,.9)(0,.25){9}{\vrule height2pt}
\put(2.8,1.2){\begin{footnotesize}$s\!-\!1$\end{footnotesize}}
\put(5.2,1.2){\begin{footnotesize}$s$\end{footnotesize}}
\multiput(-.1,3)(.25,0){45}{\hbox to 2pt{\hrulefill }}
\multiput(-.1,4.5)(.25,0){21}{\hbox to 2pt{\hrulefill }}
\multiput(-.1,5)(.25,0){17}{\hbox to 2pt{\hrulefill }}
\put(-.85,4.4){\begin{footnotesize}$y_s$\end{footnotesize}}
\put(-2.5,4.9){\begin{footnotesize}$y_s+|\la_{r}|$\end{footnotesize}}
\put(-1.6,2.9){\begin{footnotesize}$v_{r}(f)$\end{footnotesize}}
\put(.2,1.2){\begin{footnotesize}$0$\end{footnotesize}}
\put(2,.5){\vector(-1,0){2}}\put(2,.5){\vector(1,0){2}}
\put(7,.5){\vector(-1,0){2}}\put(7,.5){\vector(1,0){4}}
\put(1,-.2){\begin{footnotesize}$N_{r}^-(r_s)$\end{footnotesize}}
\put(6.2,-.2){\begin{footnotesize}$N_{r}^-(Q_s(\phi_{r})^s)$\end{footnotesize}}
\put(0,7){\line(2,-1){9}}
\put(-2.8,6.8){\begin{footnotesize}$y_s+s|\la_{r}|$\end{footnotesize}}
\put(-.1,7){\line(1,0){.2}}
\put(9,2){\begin{footnotesize}$L$\end{footnotesize}}
\end{picture}
\end{figure}

A prime ideal $\P\mid\p$ satisfies one and only one of the following conditions:\medskip

\begin{enumerate}
 \item[(i)] $\ty\mid F_\P$.
\item[(ii)] $\ty\nmid F_\P$, but $\ty'\mid F_\P$, for $\ty'=\op{Trunc}_{i-1}(\ty)$, and some maximal $1\le i< r$.
\item[(iii)] $\op{Trunc}_0( \ty)\nmid F_\P$, or equivalently, $\psi_0\nmid \overline{F}_\P$.
\end{enumerate}\medskip

We shall prove the inequality (\ref{aim}) by an independent argument in each case. We denote throughout the proof: $e=e_0\cdots e_{r-1}$.\bigskip

\noindent{\bf Case (i): $\ty\mid F_{\P}$}\medskip

By \cite[Thm. 3.1]{HN}, for some slope $\mu$ of $N_{r}^-(f)$, we have:
\begin{equation}\label{thmpol}
w_\P(\phi_{r}(\t))=(V_{r}+|\mu|)/e,
\end{equation}
and $N_r(F_\P)$ is one-sided of slope $\mu$. Consider the type $\tilde{\ty}:=(\ty;(\phi_r,\mu,\psi))$, where $\psi$ is the monic irreducible factor of $R_{\mu,r}(F_\P)$. By construction, $\tilde{\ty}\mid F_\P$.

If $|\mu|\ge|\la_{r}|$, then Proposition \ref{vgt} applied to the type $\tilde{\ty}$ and the polynomial $Q_s(x)\phi_{r}(x)^s$ shows that (see Figure \ref{fig4}):
$$
w_\P(Q_s(\t)\phi_{r}(\t)^s)\ge (y_s+s|\mu|)/e.
$$
By (\ref{thmpol}), we get $w_\P(Q_s(\t))\ge H_s$, as desired.

If $|\mu|<|\la_{r}|$, we apply Proposition \ref{vgt} to the type $\tilde{\ty}$ and the polynomial $r_s(x)$ and we get  (see Figure \ref{fig4}):
$$
w_\P(Q_s(\t)\phi_{r}(\t)^s)\stackrel{(\ref{residue})}=w_\P(r_s(\t))\ge \left(y_s+|\la_{r}|+(s-1)|\mu|\right)/e.
$$
By (\ref{thmpol}), we get in this case a stronger inequality:
$$w_\P(Q_s(\t))\ge H_s+\left(|\la_{r}|-|\mu|\right)/e.
$$

Summing up, we get in Case (i):
\begin{equation}\label{stronger1}
w_\P(Q_s(\t))\ge H_s+\max\left\{0,\left(|\la_{r}|-|\mu|\right)/e\right\}.
\end{equation}\medskip

\begin{figure}\caption{}\label{fig4}
\setlength{\unitlength}{5.mm}
\begin{picture}(11,11)
\put(5.85,3.85){$\bullet$}\put(10.85,2.85){$\bullet$}
\put(4.85,4.35){$\bullet$}\put(4.85,6.85){$\bullet$}
\put(0,0){\line(0,1){11}}\put(-1,1){\line(1,0){13}}
\put(11,3){\line(-5,1){5}}\put(11,3.03){\line(-5,1){5}}
\put(6,4){\line(-1,3){1}}\put(6,4.03){\line(-1,3){1}}
\put(5,4.5){\line(2,-1){1}}\put(5,4.53){\line(2,-1){1}}
\multiput(5,.9)(0,.25){25}{\vrule height2pt}
\multiput(11,.9)(0,.25){9}{\vrule height2pt}
\put(10.3,.35){\begin{footnotesize}$\ord_\ty(f)$\end{footnotesize}}
\put(4.9,.35){\begin{footnotesize}$s$\end{footnotesize}}
\multiput(-.1,3)(.25,0){45}{\hbox to 2pt{\hrulefill }}
\multiput(-.1,4.5)(.25,0){21}{\hbox to 2pt{\hrulefill }}
\put(-.8,4.4){\begin{footnotesize}$y_s$\end{footnotesize}}
\put(-1.6,2.9){\begin{footnotesize}$v_{r}(f)$\end{footnotesize}}
\put(.2,.35){\begin{footnotesize}$0$\end{footnotesize}}
\put(4.4,9){\begin{footnotesize}$N_{r}^-(Q_s(\phi_{r})^s)$\end{footnotesize}}
\put(0,7){\line(2,-1){9}}\put(-.1,7){\line(1,0){.2}}
\put(-2.8,6.9){\begin{footnotesize}$y_s+s|\la_{r}|$\end{footnotesize}}
\put(0,9.6){\line(1,-1){7}}\put(-.1,9.6){\line(1,0){.2}}
\put(-2.5,9.5){\begin{footnotesize}$y_s+s|\mu|$\end{footnotesize}}
\put(9,2){\begin{footnotesize}$L$\end{footnotesize}}
\put(4.8,-.6){\begin{footnotesize}$|\mu|\ge|\la_{r}|$\end{footnotesize}}
\end{picture}\qquad\qquad\qquad
\begin{picture}(5,10.5)
\put(-.15,9.85){$\bullet$}\put(1.85,5.85){$\bullet$}
\put(3.85,7.85){$\bullet$}\put(2.85,5.85){$\bullet$}
\put(0,0){\line(0,1){11}}\put(-1,1){\line(1,0){8}}
\put(2,6.03){\line(-1,2){2}}\put(2,6){\line(-1,2){2}}
\put(2,6){\line(1,0){1}}\put(2,6.03){\line(1,0){1}}
\put(3,6){\line(1,2){1}}\put(3,6.03){\line(1,2){1}}
\put(2,6){\line(2,-1){2}}\put(2,6.03){\line(2,-1){2}}
\put(3.85,4.85){$\bullet$}\put(4.75,4.35){$\times$}
\multiput(4,.9)(0,.25){28}{\vrule height2pt}
\multiput(5,.9)(0,.22){17}{\vrule height2pt}
\put(3.4,.35){\begin{footnotesize}$s\!-\!1$\end{footnotesize}}
\put(4.9,.35){\begin{footnotesize}$s$\end{footnotesize}}
\multiput(-.1,3)(.25,0){27}{\hbox to 2pt{\hrulefill }}
\multiput(-.1,4.5)(.25,0){21}{\hbox to 2pt{\hrulefill }}
\put(-.8,4.4){\begin{footnotesize}$y_s$\end{footnotesize}}
\put(-1.6,2.9){\begin{footnotesize}$v_{r}(f)$\end{footnotesize}}
\put(.2,.35){\begin{footnotesize}$0$\end{footnotesize}}
\put(2,9){\begin{footnotesize}$N_{r}^-(r_s)$\end{footnotesize}}
\put(0,6){\line(4,-1){6}}\put(-.1,6){\line(1,0){.2}}
\put(-5.5,5.9){\begin{footnotesize}$y_s+|\la_{r}|+(s\!-\!1)|\mu|$\end{footnotesize}}
\put(0,7){\line(2,-1){6}}\put(-.1,7){\line(1,0){.2}}
\put(-2.8,6.8){\begin{footnotesize}$y_s+s|\la_{r}|$\end{footnotesize}}
\put(5.6,3.6){\begin{footnotesize}$L$\end{footnotesize}}
\put(2.5,-.6){\begin{footnotesize}$|\mu|<|\la_{r}|$\end{footnotesize}}
\end{picture}
\end{figure}

\noindent{\bf Case (ii): $\ty\nmid F_{\!\P}$, $\ty'\mid F_{\P}$, for $\ty'=\mbox{\bf Trunc}_{i-1}(\ty)$, and $i$ maximal, $1\le i< r$}\medskip

As above, $N_i(F_\P)$ is one-sided of slope $\mu$, one of the slopes of $N_i^-(f)$, and
\begin{equation}\label{thmpol2}
w_\P(\phi_i(\t))=(V_i+|\mu|)/(e_0\cdots e_{i-1}).
\end{equation}
On the other hand, the arguments in the proof of \cite[Prop. 3.8]{newapp} show that
\begin{equation}\label{wphit}
w_\P(\phi_j(\t))=\dfrac{m_j}{m_i}\,\dfrac{V_i+\min\{|\la_i|,|\mu|\}}{e_0\cdots e_{i-1}},\quad i<j\le r.
\end{equation}

Since $r_s(x)=\sum_{0\le t<s}a_t(x)\phi_{r}(x)^t$, there exists $0\le t<s$ such that $w_\P(r_s(\t))\ge w_\P(a_t(\t)\phi_{r}(\t)^t)$; thus, by (\ref{residue}),
\begin{equation} \label{ak}
w_\P(Q_s(\t))=w_\P(r_s(\t))-sw_\P(\phi_{r}(\t))\ge w_\P(a_t(\t))-(s-t)w_\P(\phi_{r}(\t)). 
\end{equation}
Let us show that $w_\P(a_t(\t))$ is sufficiently large.
Consider the multiadic expansion:
\begin{equation}\label{at}
a_t(x)=\sum\nolimits_{\j\in J}b_\j(x)\Phi(x)^\j,\quad \deg b_\j<m_i,
\end{equation}
where $\Phi(x)^{\j}:=\phi_i(x)^{j_i}\cdots \phi_{r-1}(x)^{j_{r-1}}$, and $0\le j_k<e_kf_k$, for all $i\le k< r$.
Fix a multiindex $\j$ such that 
\begin{equation}\label{multi}
w_\P(a_t(\t))\ge w_\P\left(b_\j(\t)\Phi(\t)^\j\right). 
\end{equation}
Since $\ty'\mid F_\P$ and $\deg b_\j<m_i$,  \cite[Props. 2.9,2.7,(1)]{HN} show that: 
\begin{equation}\label{smalldeg}
w_\P(b_\j(\t))=v_i(b_\j)/(e_0\cdots e_{i-1})=v_{r}(b_\j)/(e_0\cdots e_{r-1}).
\end{equation}
Finally, by the convexity of the Newton polygon $N_{r}^-(f)$:
\begin{equation}\label{convex}
v_{r}\left(a_t(\phi_{r})^t\right)\ge y_s+(s-t)|\la_{r}|.
\end{equation}
If we gather (\ref{ak}), (\ref{multi}), (\ref{smalldeg}), (\ref{convex}), and we use Lemma \ref{multiadic}, we get:
\begin{align*}
 w_\P(Q_s(\t))\ge &\  w_\P(a_t(\t))-(s-t)w_\P(\phi_{r}(\t))\\
\ge &\ w_\P(b_\j(\t))+w_\P(\Phi(\t)^\j)-(s-t)w_\P(\phi_{r}(\t))\\
= &\ v_{r}(b_\j)/e+w_\P(\Phi(\t)^\j)-(s-t)w_\P(\phi_{r}(\t))\\
\ge &\ (v_{r}(a_t)-v_{r}(\Phi^\j))/e+w_\P(\Phi(\t)^\j)-(s-t)w_\P(\phi_{r}(\t))\\
\ge &\ (y_s+(s-t)|\la_{r}|-tV_{r}-v_{r}(\Phi^\j))/e+w_\P(\Phi(\t)^\j)-(s-t)w_\P(\phi_{r}(\t)).
\end{align*}
If we add and substract $sV_{r}/e$ to the last term, we get:
\begin{equation}\label{AB}
 w_\P(Q_s(\t))\ge H_s+(s-t)M+N,
\end{equation}
where
$$
M:=\dfrac{V_{r}+|\la_{r}|}{e}-w_\P(\phi_{r}(\t)),\quad N:=w_\P(\Phi(\t)^\j)-\dfrac{v_{r}(\Phi^\j)}{e}.
$$
Now, (\ref{vrphir}) and (\ref{wphit}) provide a closed formula for $M$:
\begin{align*}
M=& \ 
\dfrac{|\la_{r}|}{e}+\left(\sum_{i\le k< r}\dfrac{m_{r}}{m_k}\,\dfrac{h_k}{e_1\cdots e_k}\right)-\dfrac{m_{r}}{m_i}\,\dfrac{\min\{|\la_i|,|\mu|\}}{e_0\cdots e_{i-1}}\\
=&\ \left(\sum_{i\le k\le r}\dfrac{m_{r}}{m_k}\,\dfrac{h_k}{e_1\cdots e_k}\right)-\dfrac{m_{r}}{m_i}\,\dfrac{\min\{|\la_i|,|\mu|\}}{e_0\cdots e_{i-1}}\ge0,
\end{align*}
the last inequality because $h_i/(e_1\cdots e_{i})=|\la_i|/(e_0\cdots e_{i-1})$. Also,
(\ref{thmpol2}), (\ref{vrphis}), (\ref{wphit}) and (\ref{vrphir}) provide a closed formula for $N$:
\begin{align*}
N&=\sum\nolimits_{i\le k< r}j_k\left(w_\P(\phi_k(\t))-\dfrac{v_{r}(\phi_k)}{e}\right)\\&  =j_i\dfrac{|\mu|-|\la_i|}{e_0\cdots e_{i-1}}
 +\sum_{i< k< r}j_k\left( 
\dfrac{m_k}{m_i}\,\dfrac{V_i+\min\{|\la_i|,|\mu|\}}{e_0\cdots e_{i-1}}-\dfrac{V_k}{e_1\cdots e_{k-1}}-\dfrac{h_k}{e_1\cdots e_{k}}\right)\\
&=j_i\dfrac{|\mu|-|\la_i|}{e_0\cdots e_{i-1}}+\sum_{i< k< r}j_k\left(\dfrac{m_k}{m_i}\,\dfrac{\min\{|\la_i|,|\mu|\}}{e_0\cdots e_{i-1}}-\sum_{i\le j\le k} 
\dfrac{m_k}{m_j}\,\dfrac{h_j}{e_1\cdots e_j}\right).
\end{align*}

Since $M\ge0$ and $s>t$, we deduce from (\ref{AB}) that:
$ w_\P(Q_s(\t))\ge H_s+M+N$. Thus, we need to find lower bounds for $M+N$. 
Let us calculate first the sum of all terms of $M+N$ involving $h_i,\,|\la_i|=h_i/e_i$ and $|\mu|$. If $|\mu|\ge|\la_i|$, this partial sum is equal to $j_i(|\mu|-|\la_i|)/(e_0\cdots e_{i-1})\ge0$, whereas 
for $|\mu|<|\la_i|$ we get 
$$
\dfrac{m_{r}}{m_i}\,\dfrac{|\la_i|-|\mu|}{e_0\cdots e_{i-1}}-
\sum_{i\le k< r}j_k\dfrac{m_k}{m_i}\,\dfrac{|\la_i|-|\mu|}{e_0\cdots e_{i-1}}\ge
\dfrac{|\la_i|-|\mu|}{e_0\cdots e_{i-1}},
$$
because $(m_{r}/m_i)-\sum_{i\le k< r}j_k(m_k/m_i)\ge 1$, by Lemma \ref{combasic}.

Finally, the partial sum of the terms of $M+N$ involving $h_j$, for each $i<j\le r$, is equal to:
$$
\dfrac{m_{r}}{m_j}\,\dfrac{h_j}{e_1\cdots e_j}-
\sum_{j\le k< r}j_k\dfrac{m_k}{m_j}\,\dfrac{h_j}{e_1\cdots e_j}\ge
\dfrac{h_j}{e_1\cdots e_j},
$$
because $(m_{r}/m_j)-\sum_{j\le k< r}j_k(m_k/m_j)\ge 1$, by Lemma \ref{combasic}. Summing up, we have proven that
\begin{equation}\label{stronger2}
w_\P(Q_s(\t)) \ge \ H_s+\max\left\{0,\dfrac{|\la_i|-|\mu|}{e_0\cdots e_{i-1}}\right\}+\sum_{i<j\le r}\dfrac{h_j}{e_1\cdots e_j}>H_s.
\end{equation}

\noindent{\bf Case (iii): $\op{\bf Trunc}_0(\ty)\nmid F_\P$, or equivalently, $\psi_0\nmid \overline{F}_\P$.}\medskip

The proof is similar to the previous case, but the arguments are now simplified because $w_\P(\phi_j(\t))=0$, for all $1\le j\le r$. The formula (\ref{ak}) now gives:
$$
w_\P(Q_s(\t))=w_\P(r_s(\t))\ge w_\P(a_t(\t)),
$$
for some $0\le t<s$. If we consider the multiadic expansion (\ref{at}) for $i=1$, there exists a multiindex $\j=(j_1,\dots,j_{r-1})$ such that 
$$
w_\P(a_t(\t))\ge w_\P(b_\j(\t)\Phi(\t)^\j)=w_\P(b_\j(\t)). 
$$
Clearly, $w_\P(b_\j(\t))\ge v_1(b_\j)$. On the other hand, since $\deg b_\j<m_1$, the recursive definition of $v_2,\dots,v_{r}$ leads to
$v_1(b_\j)=v_{r}(b_\j)/(e_0\cdots e_{r})=v_{r}(b_\j)/e$.
These inequalities, together with Lemma \ref{multiadic} and (\ref{convex}), show that:
\begin{align*}
w_\P(Q_s(\t))\ge &\  v_{r}(b_\j)/e\ge \left(v_{r}(a_t)-v_{r}(\Phi^\j)\right)/e
\ge \left(y_s+(s-t)|\la_{r}|-tV_{r}-v_{r}(\Phi^\j)\right)/e\\
\ge & \ H_s+\left((s-t)(V_{r}+|\la_{r}|)-v_{r}(\Phi^\j)\right)/e
\ge H_s+\left(V_{r}+|\la_{r}|-v_{r}(\Phi^\j)\right)/e.
\end{align*}
Now, by (\ref{vrphiHN}):
$$\dfrac{V_{r}+|\la_{r}|-v_{r}(\Phi^\j)}e=
\sum_{1\le j\le r}\left(\dfrac{m_{r}}{m_j}-\sum_{j\le k< r}j_k\dfrac{m_k}{m_j}\right)\dfrac{h_j}{e_1\cdots e_j}\ge
\sum_{1\le j\le r}\dfrac{h_j}{e_1\cdots e_j},
$$
the last inequality by Lemma \ref{combasic}.
Thus, we obtain in this case:
\begin{equation}\label{stronger3}
 w_\P(Q_s(\t))\ge H_s+\sum\nolimits_{1\le j\le r}\dfrac{h_j}{e_1\cdots e_j}>H_s.
\end{equation}
\end{proof}

\subsection{Residual polynomials of quotients}\label{subsecRespolQ}
In this section we compute the resi\-dual polynomial $R_{r}(Q_s)$ of a $\phi_{r}$-quotient of $f(x)$. To this end, we recall first the general construction of the operator $R_{\lambda,i}$, of order $i>0$, with respect to a type $\ty$, of order $i-1$, a representative $\phi_i$ of $\ty$, and a negative rational number $\lambda$.

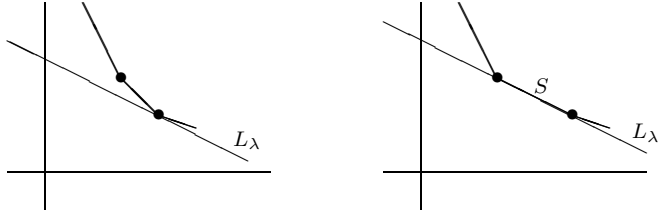
\begin{figure}\caption{$\lambda$-component of a polygon. $L_\lambda$ is the line of slope $\lambda$ having first contact with the polygon from below.}\label{figComponent}
\begin{center}
\setlength{\unitlength}{5mm}
\begin{picture}(15,5)
\put(2.85,1.45){$\bullet$}\put(1.85,2.45){$\bullet$}
\put(-1,0.1){\line(1,0){7}}\put(0,-.9){\line(0,1){5.5}}
\put(3,1.6){\line(-1,1){1}}\put(3.02,1.6){\line(-1,1){1}}
\put(3,1.6){\line(3,-1){1}}\put(3.02,1.6){\line(3,-1){1}}
\put(2,2.6){\line(-1,2){1}}\put(2.02,2.6){\line(-1,2){1}}
\put(5.4,.4){\line(-2,1){6.4}}\put(5,.8){\begin{footnotesize}$L_{\lambda}$\end{footnotesize}}
\put(13.85,1.45){$\bullet$}\put(11.85,2.45){$\bullet$}
\put(9,0.1){\line(1,0){7}}\put(10,-.9){\line(0,1){5.5}}
\put(14,1.6){\line(-2,1){2}}\put(14.02,1.6){\line(-2,1){2}}
\put(14,1.6){\line(3,-1){1}}\put(14.02,1.6){\line(3,-1){1}}
\put(12,2.6){\line(-1,2){1}}\put(12.02,2.6){\line(-1,2){1}}
\put(16,.6){\line(-2,1){7}}\put(15.6,1){\begin{footnotesize}$L_{\lambda}$\end{footnotesize}}
\put(13,2.2){\begin{footnotesize}$S$\end{footnotesize}}
\end{picture}
\end{center}
\end{figure}

\begin{definition}\label{lacomponent}
Let $\lambda\in\Q_{<0}$ and $N$ a Newton polygon. We define the \emph{$\lambda$-compo\-nent} of $N$ to be $S_{\lambda}(N):=\{(x,y)\in N\mid y-\lambda x\mbox{ is minimal}\}$. 
If $N$ has a side $S$ of slope $\lambda$, then $S_{\lambda}(N)=S$; otherwise, $S_\lambda(N)$ is a vertex of $N$ (see Figure \ref{figComponent}).  
\end{definition}

Let $\lambda=-h/e$, with $h,e$ positive coprime integers. 
Let $g(x)=\sum_{0\le s}a_s(x)\phi_i(x)^s$ be the canonical $\phi_i$-expansion of a non-zero polynomial $g(x)\in \oo_\p[x]$. Let $S$ be the $\lambda$-component of $N_i(g)$, and let $s_0$ be the abscissa of the left end point of $S$.  The points of integer coordinates lying on $S$ have abscissa  $s_j:=s_0+je$, for $0\le j\le d$, where $d:=d(S)$ is the \emph{degree} of $S$. 

For each $0\le s\le \ell(N_i^-(g))$, consider the residual coefficient $c_{s}\in\F_i$, defined as:
\begin{equation}\label{rescoeff}
c_{s}:=\begin{cases}
0,&\mbox{ if $(s,v_i(a_{s}\phi^{s}))$ lies above }N_i^-(g),\\
z_{i-1}^{t_{i-1}(s)}R_{i-1}(a_{s})(z_{i-1}),&\mbox{ if $(s,v_i(a_{s}\phi^{s}))$ lies on }N_i^-(g),             
\end{cases}
\end{equation}
where $t_0(s):=0$, and $t_{i-1}(s)$ is described in \cite[Def. 2.19]{HN} for $i>1$. 

\begin{definition}\label{respol}
The \emph{residual polynomial} of $g$ with respect to $(\ty,\phi_i,\lambda)$ is:
$$
R_{\lambda,i}(g)(y):=c_{s_0}+c_{s_1}y+\cdots+c_{s_d}y^d\in \F_i[y].
$$
Since $c_{s_0}c_{s_d}\ne0$, this polynomial has degree $d=d(S)$ and it is not divisible by $y$.
\end{definition}

We now go back to the situation described in Theorem \ref{denquot}. Recall that $S$ is the $\lambda_r$-component of $N_r^-(f)$, and $1\le s\le \ell(N_r^-(f))$ is an abscissa belonging to the projection of $S$ to the horizontal axis. Denote by $d'=d(S)$ the degree of $S$, and let $(b,u)$ be the right end point of $S$.
By Definition \ref{respol},
$$
R_{r}(f)(y)=c_{b-d'e_r}+c_{b-(d'-1)e_r}y+\cdots +c_{b}y^{d'},
$$ 
The points of $S$ having integer coordinates are marked in Figure \ref{figResQ} with $\circ$; among them, those belonging to the cloud of points $(t,v_r(a_t\phi_r^t))$ are marked with $\bullet$. 

\begin{lemma}\label{coincidence}
Let $N=N_{r}^-(f)$, $N'=N_{r}^-(Q_{s}(\phi_{r})^{s})$.
Let $d$ be the greatest integer $0\le d\le \lfloor (b-s)/e_r\rfloor$ such that $c_{b-de_r}\ne0$. Then,
the $\la_{r}$-component of $N'$ has left end point $(b-de_r,u+dh_{r})$ and right end point $(b,u)$. Moreover, 
$$
R_{r}(Q_{s}(\phi_{r})^{s})(y)=c_{b-de_r}+c_{b-(d-1)e_r}y+\cdots +c_{b}y^{d}.
$$ 
\end{lemma}

\begin{proof}
Figure \ref{figResQ2} shows the shape of $N'=N_{r}^-(Q_s(\phi_{r})^{s})$. By definition,
$$
Q_{s}(x)\phi_{r}(x)^{s}=\sum\nolimits_{s\le t}a_t(x)\phi_{r}(x)^{t},
$$
and $N'$ is the lower convex hull of the cloud of points of the plane with coordinates $(t,v_r(a_t\phi_r^t))$, for $s\le t$. Clearly, 
$N'\cap \left([b-de_r,\infty)\times\R\right)=N\cap \left([b-de_r,\infty)\times\R\right)$, and
the residual coefficients of $N'$ are $c'_t=c_{t}$, for all integer abscissas $s\le t\le \ell(N')=\ell(N)$.
The lemma is an immediate consequence of this fact.
\end{proof}

\begin{figure}\caption{}\label{figResQ}
\setlength{\unitlength}{5.mm}
\begin{picture}(14,8)
\put(1.85,5.85){$\bullet$}\put(10.85,2.85){$\bullet$}
\put(3.35,5.35){$\bullet$}\put(4.85,4.85){$\circ$}\put(6.35,4.35){$\circ$}
\put(7.85,3.85){$\bullet$}\put(9.35,3.35){$\circ$}
\put(0,0){\line(0,1){8}}\put(-1,1){\line(1,0){15}}
\put(2,6.03){\line(-1,2){1}}\put(2,6){\line(-1,2){1}}
\put(2,6){\line(3,-1){9}}\put(2,6.03){\line(3,-1){9}}
\put(11,3){\line(5,-1){2}}\put(11,3.03){\line(5,-1){2}}
\put(5.75,4.5){$\times$}
\put(6,6.5){\begin{footnotesize}$N_{r}^-(f)$\end{footnotesize}}
\put(7,4.7){\begin{footnotesize}$\la_{r}$\end{footnotesize}}
\multiput(2,.9)(0,.25){21}{\vrule height2pt}
\multiput(6,.9)(0,.25){16}{\vrule height2pt}
\multiput(8,.9)(0,.25){13}{\vrule height2pt}
\multiput(11,.9)(0,.25){9}{\vrule height2pt}
\put(1.4,.3){\begin{footnotesize}$b\!-\!d'e_r$\end{footnotesize}}
\put(7.5,.3){\begin{footnotesize}$b\!-\!de_r$\end{footnotesize}}
\put(5.8,.3){\begin{footnotesize}$s$\end{footnotesize}}
\put(10.9,.3){\begin{footnotesize}$b$\end{footnotesize}}
\put(-.6,2.9){\begin{footnotesize}$u$\end{footnotesize}}
\put(-.5,.4){\begin{footnotesize}$0$\end{footnotesize}}
\put(-.8,4.6){\begin{footnotesize}$y_{s}$\end{footnotesize}}
\multiput(-.1,4.7)(.25,0){25}{\hbox to 2pt{\hrulefill }}
\multiput(-.1,3)(.25,0){45}{\hbox to 2pt{\hrulefill }}
\end{picture}
\end{figure}
\begin{figure}\caption{}\label{figResQ2}
\setlength{\unitlength}{5.mm}
\begin{picture}(14,8)
\put(6.3,5.35){$\bullet$}\put(5.8,6.7){$\bullet$}
\put(8,4){\line(-1,1){1.6}}\put(8,4.03){\line(-1,1){1.6}}
\put(6.5,5.4){\line(-1,3){.45}}\put(6.5,5.42){\line(-1,3){.45}}
\put(10.85,2.85){$\bullet$}
\put(6.35,4.35){$\circ$}
\put(7.85,3.85){$\bullet$}\put(9.35,3.35){$\circ$}
\put(0,0){\line(0,1){7.5}}\put(-1,1){\line(1,0){15}}
\put(11,3){\line(-3,1){5}}\put(11,3.03){\line(-3,1){3}}
\put(11,3){\line(5,-1){2}}\put(11,3.03){\line(5,-1){2}}
\put(5.7,4.5){$\times$}
\put(8,6){\begin{footnotesize}$N_{r}^-(Q_{s}(\phi_{r})^{s})$
\end{footnotesize}}
\put(8.8,4){\begin{footnotesize}$\la_{r}$\end{footnotesize}}
\multiput(5.95,.9)(0,.25){24}{\vrule height2pt}
\multiput(8,.9)(0,.25){13}{\vrule height2pt}
\multiput(11,.9)(0,.25){9}{\vrule height2pt}
\put(7.4,.3){\begin{footnotesize}$b\!-\!de_r$\end{footnotesize}}
\put(5.8,.3){\begin{footnotesize}$s$\end{footnotesize}}
\put(10.9,.3){\begin{footnotesize}$b$\end{footnotesize}}
\put(-.6,2.9){\begin{footnotesize}$u$\end{footnotesize}}
\put(-.5,.4){\begin{footnotesize}$0$\end{footnotesize}}
\put(-.8,4.6){\begin{footnotesize}$y_{s}$\end{footnotesize}}
\multiput(-.1,4.7)(.25,0){25}{\hbox to 2pt{\hrulefill }}
\multiput(-.1,3)(.25,0){45}{\hbox to 2pt{\hrulefill }}
\end{picture}
\end{figure}

\begin{corollary}\label{RQ}
If we convene that $\ell_0:=0$, then,
$$
R_{r}(Q_{s})(y)=
(z_{r-1})^{s\ell_{r-1}V_{r}/e_{r-1}}\left(c_{b-de_r}+c_{b-(d-1)e_r}y+\cdots +c_{b}y^{d}\right).
$$ 
\end{corollary}

\begin{proof}
By \cite[Thm. 2.26]{HN}, we have $R_{r}(Q_{s})=R_{r}(Q_{s}(\phi_{r})^{s})/R_{r}(\phi_{r})^{s}$
in $\F_{r}[y]$. On the other hand, $R_{r}(\phi_{r})=(z_{r-1})^{-\ell_{r-1}V_{r}/e_{r-1}}\in\F_{r}^*$ is a non-zero constant polynomial, which was calculated along the proof of \cite[Thm. 2.11]{HN}.
Thus, the corollary follows from Lemma \ref{coincidence}.
\end{proof}

\begin{corollary}\label{applications}
With the above notation, let $\psi_r\in\F_r[y]$ be a monic irreducible factor of $R_r(f)$, and consider the type $\ty':=(\ty;(\phi_r,\lambda_r,\psi_r))$. Let $\P\mid\p$ be a prime ideal of $B$ such that $\ty'\mid F_\P$, and suppose $0\le b-s<e_rf_r$. Then, $w_\P(Q_s(\t))=H_s$.
\end{corollary}

\begin{proof}
By Corollary \ref{RQ}, $\deg R_r(Q_s)\le (b-s)/e_r<f_r$; thus, $\psi_r\nmid R_r(Q_s)$. Hence,  $\ty'\nmid Q_s(x)$, and Proposition \ref{vgt} shows that (see Figure \ref{figResQ2})
$$
w_\P(Q_s(\t)\phi_r(\t)^s)=(y_s+s|\lambda_r|)/(e_0\cdots e_{r-1}).
$$
On the other hand, $w_\P(\phi_r(\t))=(V_r+|\lambda_r|)/(e_0\cdots e_{r-1})$, by Corollary \ref{previous}. Therefore, $w_\P(Q_s(\t))=w_\P(Q_s(\t)\phi_r(\t)^s)-sw_\P(\phi_r(\t))=H_s$. 
\end{proof}

\section{Quotients and local integral bases}\label{secIB}
Consider again the local context of section \ref{secOkutsu}. Let $k$ be a local field, and let $v,\,\oo,\,\m,\,\pi,\,\F$, be as in that section. Also, let $F(x)\in\oo[x]$ be a monic irreducible separable polynomial, and let $\t,\, L,\,\oo_L,\,\m_L,\,\F_L$, be as in section \ref{secOkutsu}.  
We indicate with a bar, $\raise.8ex\hbox{---}\colon \oo_L\longrightarrow \F_L$, the reduction modulo $\m_L$ homomorphism. Denote: 
$$e:=e(L/k), \quad f:=f(L/k),\quad n_L:=[L\colon k]=\deg F=ef.$$

\subsection{Local bases in standard form}
Let $U_L$ be the group of units of $\oo_L$, and 
$$
\B:=\B_L:=\{\alpha\in \oo_L\mid 0\le v(\alpha)<1\}.
$$

\begin{lemma}\label{standard}
Let $\bb\subset \B$ be a finite subset. Split $\bb$ into the disjoint union:
$$
\bb=\bigcup\nolimits_{0\le t<e}\bb_t,\qquad \bb_t:=\{\alpha\in \bb\mid v(\alpha)=t/e\},
$$
and suppose that the following two conditions hold:
\begin{enumerate}
\item[(a)] $\#\bb_t=f$,
\item[(b)] for some $\omega\in \bb_t$, the family $\,\overline{\om^{-1}\bb_t}\,$ is an $\F$-basis of $\F_L$, 
\end{enumerate}
for all $0\le t<e$. Then, $\bb$ is an $\oo$-basis of $\oo_L$.
\end{lemma}

\begin{proof}
Let $M\subset \oo_L$ be the $\oo$-module generated by the elements in $\bb$. By condition (a), $\#\bb=n_L$.
By condition (b), $\bb\otimes_\oo\F$ is an $\F$-linearly independent family, hence an $\F$-basis, of $\oo_L\otimes_\oo\F$. Therefore, $M=\oo_L$ by Nakayama's lemma. 
\end{proof}

Condition (b) of the lemma does not depend of the choice of the element $\om\in\bb_t$. Actually, we can replace $\om$ by any element in $\oo_L$ having valuation $t/e$. 

\begin{definition}
An $\oo$-basis $\bb$ of $\oo_L$ is said to be in \emph{standard form} if $\bb\subset \B$ and, the two conditions of Lemma \ref{standard} are satisfied, for all $0\le t<e$.
\end{definition}

\noindent{\bf Examples. }
\begin{enumerate}
\item Suppose $U\subset U_L$ is a family of units such that $\overline{U}$ is an $\F$-basis of $\F_L$. Take $\om_0,\dots,\om_{e-1}\in\oo_L$ such that
$v(\om_t)=t/e$, for all $0\le t<e$. Then, $\bb:=\cup_{0\le t<e}\om_tU$ \,is an $\oo$-basis of $\oo_L$ in standard form.
\item The Okutsu basis $\bb=\{g_m(\t)/\pi^{\mu_m}\mid 0\le m<n_L\}$, described at the end of section \ref{secOkutsu} is in standard form \cite[I,Prop. 2]{Ok}, \cite[Prop. 4.28]{HN}.
\end{enumerate}

\begin{definition}
We define a \emph{ star operation} between elements of $\oo_L\setminus\{0\}$, by:
$$
\alpha\star\beta:=\alpha\beta/\pi^{\lfloor v(\alpha\beta)\rfloor}\in\B.
$$ 
It is clearly associative and commutative.
\end{definition}


\begin{lemma}\label{multiplybasis}
Let $\bb$ be an $\oo$-basis of $\oo_L$ in standard form. For any $\om\in \oo_L\setminus\{0\}$, the set 
$\,\om\star\bb:=\{\om\star\alpha\mid\alpha\in \bb\}$ is an $\oo$-basis of $\oo_L$ in standard form.
\end{lemma}

\begin{proof}
Let $\Pi\in\oo_L$ be a uniformizer, and write $\om=\Pi^m\eta$, for some unit $\eta\in U_L$  and some exponent $0\le m$. Clearly, $\om\star\alpha=\Pi\star\cdots\star\Pi\star\eta\star\alpha$, for all $\alpha\in\bb$. Thus, it is sufficient to check that $\eta\star\bb$ and $\Pi\star\bb$ are bases in standard form. For $\eta\star\bb=\eta\bb$ this is obvious. For $\Pi\star\bb$ we have
$$
\Pi\star\bb=\bigcup\nolimits_{0\le t<e}\bb'_t, \quad \bb'_t:=\left\{\begin{array}{ll}
\Pi\bb_{t-1},&\mbox{ if }t>0,\\
(\Pi/\pi)\bb_{e-1},&\mbox{ if }t=0.
\end{array}
\right.
$$
Clearly, the sets $\bb'_t$ satisfy the conditions of Lemma \ref{standard}, for all $0\le t<e$. 
\end{proof}

\subsection{Quotients and local bases}
Let $\ty=(\psi_0;(\phi_1,\lambda_1,\psi_1);\cdots,(\phi_r,\lambda_r,\psi_r))$ be an $F$-complete type of order $r$. By \cite[Cor. 3.8]{HN}, we have $e=e_0e_1\cdots e_r$, $f=f_0f_1\cdots f_r$.

To this type $\ty$ we may attach several rational functions in $k(x)$ \cite[Sec. 2.4]{HN}.
Let $\,\pi_0(x)=1$, $\pi_1(x)=\pi$. We define recursively for all $1\le i\le r$,
\begin{equation}\label{ratfracs}
\Phi_i(x)=\dfrac{\phi_i(x)}{\pi_{i-1}(x)^{V_i/e_{i-1}}},\qquad
\gamma_i(x)=\dfrac{\Phi_i(x)^{e_i}}{\pi_i(x)^{h_i}},\qquad
\pi_{i+1}(x)=\dfrac{\Phi_i(x)^{\ell_i}}{\pi_i(x)^{\ell'_i}}.
\end{equation}
These rational functions can be written as a product of powers of $\pi,\phi_1(x),\dots,\phi_r(x)$, with integer exponents. Recall that $\ell_i$, $\ell'_i$ are integers satisfying the identity $\ell_ih_i-\ell'_ie_i=1$ (section \ref{secOkutsu}).

The type $\ty$ determines as well a chain of extensions of the residue field of $k$:
$$
\F=\F_0\subset\F_1\subset\cdots\subset\F_{r+1},\qquad \F_{i+1}=\F[z_0,\dots,z_i], \ 0\le i\le r.
$$
The residue field $\F_L$ can be identified to the field $\F_{r+1}$.
More precisely, in \cite[(27)]{HN} we construct an explicit isomorphism
\begin{equation}\label{gamma}
\gamma\colon \F_{r+1}\lra \F_L, \qquad z_0\mapsto \overline{\t}, \ z_1\mapsto \gb1,\ \dots,\ z_r\mapsto \gb{r},
\end{equation}
where $\gamma_i(x)\in k(x)$ are the rational functions defined in (\ref{ratfracs}). 

We denote by $\red_L \colon\oo_L\lra \F_{r+1}$,  the reduction map obtained by composition of the canonical reduction map with the inverse of this isomorphism:
$$
\red_L\colon \oo_L\lra \F_L \stackrel{\gamma^{-1}}\lra \F_{r+1}.
$$

The following proposition is easily deduced from \cite[Prop. 3.5]{HN}.

\begin{proposition}\label{red}Let $F,\,\t,\,L,\,\ty,$ be as above.
Let $g(x)\in\oo[x]$ be a non-zero polynomial, and let $(s,u)$
be the left end point of the $\lambda_r$-component of $N_r(g)$ (Definition \ref{lacomponent}). If $\ty\nmid g$, we have 
$v(g(\t))=v(\Phi_r(\t)^s\pi_r(\t)^u)$, and 
$$
\red_L\left(g(\t)/(\Phi_r(\t)^s\pi_r(\t)^u)\right)=R_r(g)(z_r)\ne0.
$$
\end{proposition}

\begin{notation}\label{not}
Let $f(x)\in\oo[x]$ be a monic separable polynomial, divisible by $F(x)$ in $\oo[x]$. 
For all $1\le i\le r$, we denote

\begin{tabular}{cl}
$(b_i,u_{i})$ &the right end point of the side of slope $\la_i$ of $N_i^-(f)$\\
$u'_i$ &$=u_{i}-b_iV_i$ \\
$Q_{i,j}(x)$ &the $(b_i-j)$-th  $\phi_i$-quotient of $f(x)$, for $0\le j<b_i$\\
$y_{i,j}$ &the ordinate of the point of $N_i^-(f)$ with abscissa $b_i-j$\\
$H_{i,j}$&$=(y_{i,j}-(b_i-j)V_i)/(e_0\cdots e_{i-1})=(u'_i+j(V_i+|\la_i|))/(e_0\cdots e_{i-1})$
\end{tabular}
\end{notation}

We emphasize that $j$ is the distance from the relevant abscissa $b_i-j$ of the quotient, to the abscissa $b_i$ of the right end point of the relevant side of $N_i^-(f)$. Figure \ref{figNewQ} illustrates the situation.

\begin{figure}\caption{}\label{figNewQ}
\setlength{\unitlength}{5.mm}
\begin{picture}(20,7.8)
\put(1.85,5.35){$\bullet$}\put(5.85,3.35){$\bullet$}
\put(0,-.5){\line(0,1){8.2}}\put(-1,.5){\line(1,0){10}}
\put(2,5.53){\line(-1,2){1}}\put(2,5.5){\line(-1,2){1}}
\put(6,3.5){\line(-2,1){4}}\put(6,3.53){\line(-2,1){4}}
\put(6,3.5){\line(5,-1){2}}\put(6,3.53){\line(5,-1){2}}
\put(3.8,4.35){$\times$}
\put(5.3,6.1){\begin{footnotesize}$N_i^-(f)$\end{footnotesize}}
\put(5,4.2){\begin{footnotesize}$\lambda_i$\end{footnotesize}}
\multiput(4,.4)(0,.25){17}{\vrule height2pt}
\multiput(6,.4)(0,.25){13}{\vrule height2pt}
\put(3.3,-.2){\begin{footnotesize}$b_i\!-\!j$\end{footnotesize}}
\put(5.9,-.2){\begin{footnotesize}$b_i$\end{footnotesize}}
\put(-.9,3.4){\begin{footnotesize}$u_{i}$\end{footnotesize}}
\put(-.5,-.1){\begin{footnotesize}$0$\end{footnotesize}}
\put(-1.1,4.5){\begin{footnotesize}$y_{i,j}$\end{footnotesize}}
\put(4.9,1.2){\begin{footnotesize}$j$\end{footnotesize}}
\put(5,1){\vector(-1,0){1}}\put(5,1){\vector(1,0){1}}
\multiput(-.1,4.55)(.25,0){17}{\hbox to 2pt{\hrulefill }}
\multiput(-.1,3.5)(.25,0){25}{\hbox to 2pt{\hrulefill }}
\put(17.35,1.8){$\bullet$}
\put(16.35,2.3){$\bullet$}\put(15.35,5.2){$\bullet$}
\put(15.5,-.5){\line(0,1){8.2}}\put(14.5,.5){\line(1,0){7}}
\put(15.5,5.5){\line(1,-3){1}}\put(15.5,5.53){\line(1,-3){1}}
\put(17.5,2){\line(-2,1){2}}\put(17.5,2.03){\line(-2,1){1}}
\put(17.5,2){\line(5,-1){2}}\put(17.5,2.03){\line(5,-1){2}}
\put(15.25,2.85){$\times$}
\put(17,6.1){\begin{footnotesize}$N_i^-(Q_{i,j})$\end{footnotesize}}
\put(16.9,2.4){\begin{footnotesize}$\lambda_i$\end{footnotesize}}
\multiput(17.5,.4)(0,.25){7}{\vrule height2pt}
\put(17.4,-.2){\begin{footnotesize}$j$\end{footnotesize}}
\put(11.7,1.8){\begin{footnotesize}$u_{i}\!-\!(b_i\!-\!j)V_i$\end{footnotesize}}
\put(15,-.1){\begin{footnotesize}$0$\end{footnotesize}}
\put(11.4,2.9){\begin{footnotesize}$y_{i,j}\!-\!(b_i\!-\!j)V_i$\end{footnotesize}}
\multiput(15.4,2)(.25,0){8}{\hbox to 2pt{\hrulefill }}
\end{picture}
\end{figure}

The aim of this section is to prove the following result.

\begin{thm}\label{Pbasis}
Let $J=\{(j_0,\dots,j_r)\in\N^{r+1}\mid 0\le j_i<e_if_i,\ \forall\,0\le i\le r\}$, and for each multiindex $\j\in J$, consider:
$$
Q_\j:=\t^{j_0} Q_{1,j_1}(\t)\star\cdots\star Q_{r,j_r}(\t)=\dfrac{\t^{j_0}Q_{1,j_1}(\t)\cdots Q_{r,j_r}(\t)}{\pi^{\lfloor H_{1,j_1}+\cdots +H_{r,j_r}\rfloor}}\in\B.
$$
Then, the family $\bb:=\{Q_\j\mid \j\in J\}$ is an $\oo$-basis of $\oo_L$ in standard form.
\end{thm}

By (\ref{frame}), $v(\t^{j_0})=0$, because $0\le j_0<e_0f_0=m_1$. Also, $v(Q_{i,j_i}(\t))=H_{i,j_i}$, for all $i$, by Corollary \ref{applications}. Thus, $Q_\j$ indeed belongs to $\B$.
In order to prove Theorem \ref{Pbasis} we need to check that the sets, $\bb_t=\{\alpha\in\bb \mid v(\alpha)=t/e\}$, satisfy the conditions of Lemma \ref{standard}. This will be shown by a recursive argument. 

Consider the filtration of \ $\B$ determined by the subsets:
$$
U_L=\B_0\subset\B_1\subset\cdots\subset \B_{r}=\B,\quad \B_i:=\{\alpha\in\B\mid e_0\cdots e_i\,v(\alpha)\in\Z\}.
$$
For $0\le i\le r$, $\B_i$ splits as the disjoint union:
$$\B_i=\bigcup\nolimits_{0\le t< e_0\cdots e_i}\B_{i,t},\qquad \B_{i,t}:=\left\{\alpha\in\B_i \mid  v(\alpha)=t/e_0\cdots e_i\right\}.$$

\begin{definition}\label{partialstandard}
We say that $B\subset \B_i$ is a \emph{level $i$ basis in standard form} if  for all $0\le t<e_0\cdots e_i$, the following two conditions are satisfied:
\begin{enumerate}
\item $\#B_t=f_0\cdots f_i$, where $B_{t}:= B\cap \B_{i,t}$,
\item For any $\om\in\B_{i,t}$, the family $\red_L(\om^{-1}B_t)$ is an $\F$-basis of $\F_{i+1}$.
\end{enumerate}
\end{definition}

\begin{lemma}\label{valuett}
Let $B\subset \B_{i-1}$ be a level $i-1$ basis in standard form, for some $1 \le i\le r$. For each $0\le t<e_0\cdots e_i$, take $0\le q_t<e_i$ such that: $q_th_i\equiv t\md{e_i}$. \vskip1.5mm

(a)  Let $\left\{\om_j\right\}_{0\le j<e_if_i}\subset\B_i$ such that $\ e_0\cdots e_i\,v(\om_j)\equiv jh_i\md{e_i}
$, for all $j$. Then, 
\begin{enumerate}
\item[(i)] If $j\not\equiv q_t \md{e_i}$, then $(\om_j\star B)\cap \B_{i,t}=\emptyset$.
\item[(ii)] For each $j=q_t+ke_i$, there exists a unique $0\le t_k<e_0\cdots e_{i-1}$, depending on $i,t,k$ and $v(\omega_j)$, such that $\om_j\star B_{t_k}\subset \B_{i,t}$.
\end{enumerate}\vskip1.5mm

(b) \ Suppose moreover that for some $\Pi\in \B_{i,1}$, $\Pi_0\in \B_{i-1,1}$, the family
$$
\epsilon_k:=\red_L\left(\Pi^{-t}\om_{q_t+ke_i}\star (\Pi_0)^{t_k}\right),\quad 0\le k<f_i,
$$
is an $\F_i$-basis of $\F_{i+1}$, for all $0\le t<e_0\cdots e_i$. Then, 
$B':=\bigcup\nolimits_{0\le j<e_if_i}\om_j\star B\subset\B_i$ 
is a level $i$ basis in standard form.
\end{lemma}

\begin{proof}Let $e_0\cdots e_i\,v(\om_j)=jh_i+p_je_i$, for some integer $p_j$. 
For any $\,0\le t'<e_0\cdots e_{i-1}$, the condition $\om_j\star B_{t'}\subset \B_{i,t}$ is equivalent to:
$$
v(\om_j)+t'/(e_0\cdots e_{i-1})\equiv t/(e_0\cdots e_i)\md{\Z},
$$
or, equivalently,
\begin{equation}\label{modei}
jh_i+p_je_i+t'e_i\equiv t \md{e_0\cdots e_i}.
\end{equation}
The condition $jh_i\equiv t\md{e_i}$, or equivalently, $j\equiv q_t\md{e_i}$, is necessary. On the other hand, if we denote $n_{i,t}:=(q_th_i-t)/e_i$, and we take $j=q_t+ke_i$, then there is a unique $0\le t'<e_0\cdots e_{i-1}$ for which (\ref{modei}) holds:
$$
t'\equiv -(n_{i,t}+kh_i+p_j)\md{e_0\cdots e_{i-1}}.
$$
This proves items (i) and (ii) of the lemma. 
Therefore, the set $B'$ splits as:
$$
B'=\bigcup\nolimits_{0\le t<e_0\cdots e_i}B'_t,\qquad
B'_t=\bigcup\nolimits_{0\le k<f_i}\om_{q_t+ke_i}\star B_{t_k}.
$$
In particular, $\#B'_t=f_i\#B_{t_k}=f_0\cdots f_i$.

Finally, by the hypothesis on $B$, for any given $t$ as above, the sets $B_{t_k}$ can be expressed as $B_{t_k}=(\Pi_0)^{t_k}U_k$, for a set $U_k\subset U_L$ such that $\red_L(U_k)$ is an $\F$-basis of $\F_i$. Hence, if the family $(\epsilon_k)_{0\le k<f_i}$ is an $\F_i$-basis of $\F_{i+1}$, then, the family
$$
\red_L(\Pi^{-t}B'_t)=\bigcup\nolimits_{0\le k<f_i}\epsilon_k\red_L(U_k)
$$
is an $\F$-basis of $\F_{i+1}$.
\end{proof}

Lemma \ref{valuett} can be applied to construct different local integral bases in standard form. Starting with $B_0=\{1,\,\t,\cdots,\t^{f_0-1}\}$, we recursively construct, for $1\le i\le r$,
$$
B_i=\bigcup\nolimits_{0\le j<e_if_i}\om_{i,j}\star B_{i-1},
$$
with $\om_{i,j}$ satisfying the conditions of Lemma \ref{valuett}. For instance, we can take $\om_{i,j}=\phi_i(\t)^j$, and we reobtain Okutsu's basis, as described in Theorem \ref{basis}.

Let us apply this idea to the quotients. For all $0\le i\le r$, consider the set 
$$
B_i:=\{\t^{j_0} Q_{1,j_1}\star\cdots\star Q_{i,j_i}\mid 0\le j_k<e_kf_k,\ \forall\,0\le k\le i\}\subset \B_i.
$$
Since $\bb=B_r$, and $\F_{r+1}\simeq\F_L$, Theorem \ref{Pbasis} is a consequence of the following result.

\begin{proposition}\label{PbasisFinal}
For all $0\le i\le r$, the set $B_i$ is a level $i$ basis in standard form.
\end{proposition}
 
\begin{proof}
We prove the proposition by induction on $i$. For $i=0$, we have
$B_0=\{1,\,\t,\,\dots,\,\t^{f_0-1}\}$, and the statement is clear.
Suppose $i>0$ and $B_{i-1}$ is a level $i-1$ basis in standard form. In order to show that $B_i$ is a level $i$ basis in standard form, we need only to check that the family $\om_j:=Q_{i,j}(\t)$, for $0\le j<e_if_i$, satisfies the conditions (a) and (b) of Lemma \ref{valuett}.

Condition (a) on $v(\om_j)$ is clearly satisfied:
$$
e_0\cdots e_i\,v(\om_j)=e_0\cdots e_i\,H_{i,j}=u'_ie_i+j(e_iV_i+h_i)\equiv jh_i\md{e_i}.
$$

Let us prove condition (b). The rational functions $\pi_i(x)$, $\pi_{i+1}(x)$ defined in (\ref{ratfracs}) satisfy \cite[Cor. 3.2]{HN}:
$$v(\pi_i(\t))=1/(e_0\cdots e_{i-1}), \quad v(\pi_{i+1}(\t))=1/(e_0\cdots e_{i}).
$$
By taking, $\Pi:=\pi_{i+1}(\t)$, $\Pi_0:=\pi_{i}(\t)$, we need only to show that, for any $0\le t<e_0\cdots e_i$, the family $(\epsilon_k)_{0\le k<f_i}$, considered in (b) is an $\F_i$-basis of $\F_{i+1}$. 

Let $0\le t<e_0\cdots e_i$, and consider the integer $0\le q_t<e_i$ of Lemma \ref{valuett}. Since $\ell_ih_i-\ell'_ie_i=1$, the integer $N:=(q_t-\ell_it)/e_i$ depends only on $i$ and $t$. Clearly, $(q_th_i-t)/e_{i}=Nh_i+\ell'_it$.

We fix an integer $0\le k<f_i$, and we let $j:=q_t+ke_{i}$. By (a) of Lemma \ref{valuett},  there is a unique $0\le t_k<e_0\cdots e_{i-1}$ such that  
$$
n_k:=H_{i,j}+t_k/(e_0\cdots e_{i-1})-t/(e_0\cdots e_{i})
$$
is a non-negative integer that depends on $i$, $t$ and $k$. We can express:
\begin{equation}\label{nk}
\begin{aligned}
e_0\cdots e_{i-1}n_k=&\ u'_i+j(V_i+|\la_i|)+
t_k-t/e_{i}\\
=&\ u'_i+jV_i+kh_i+t_k+
(q_th_i-t)/e_{i}\\
=&\ u'_i+jV_i+kh_i+t_k+Nh_i+\ell'_it.
\end{aligned}
\end{equation}

By definition, $\epsilon_k\in\F_{i+1}$, is the image under $\red_L$ of the unit 
\begin{equation}\label{unit}
\pi_{i+1}(\t)^{-t}Q_{i,j}(\t)\star\pi_i(\t)^{t_k}=
\dfrac{Q_{i,j}(\t)\pi_i(\t)^{t_k}}{\pi^{n_k}\pi_{i+1}(\t)^t}.
\end{equation}

Figure \ref{figResQ2} shows the shape of $N_i^-(Q_{i,j}(\phi_i)^{b_i-j})$. Let $d_k$ be the degree of $R_i(Q_{i,j})$. By Lemma \ref{coincidence}, $d_k\le \lfloor j/e_i\rfloor=k<f_i$, and the left end point of the $\la_i$-component of this polygon is $(b_i-d_ke_i,u_{i}+d_kh_i)$. Clearly, the Newton polygon $N_i^-(Q_{i,j})$ is the image of the former polygon under the following transformation of the plane: $$(x,y)\mapsto(x-(b_i-j),y-(b_i-j)V_i).$$ Hence, $R_i(Q_{i,j})(y)$ is not divisible by $\psi_i(y)$, and the left end point of the $\la_i$-compo\-nent of $N_i^-(Q_{i,j})$ has coordinates $(j-d_ke_i,u'_i+jV_i+d_kh_i)$. By Proposition \ref{red},
\begin{equation}\label{fraction}
\red_L\left(\dfrac{Q_{i,j}(\t)}{\Phi_i(\t)^{s}\pi_i(\t)^{u}}\right)=R_i(Q_{i,j})(z_i)\in\F_{i+1}^*,
\end{equation}
where $s:=j-d_ke_i$ and $u:=u'_i+jV_i+d_kh_i$.

We can express the unit (\ref{unit}) as the product of two units:
\begin{equation}\label{twounits}
\dfrac{Q_{i,j}(\t)\pi_i(\t)^{t_k}}{\pi^{n_k}\pi_{i+1}(\t)^t}=
\dfrac{Q_{i,j}(\t)}{\Phi_i(\t)^{s}\pi_i(\t)^{u}}\cdot \dfrac{\Phi_i(\t)^{s}\pi_i(\t)^{u+t_k}}{\pi^{n_k}\pi_{i+1}(\t)^t},
\end{equation}
and the residue class of the first factor is computed in (\ref{fraction}). If we use the following identities from (\ref{ratfracs}): 
$$
\pi_{i+1}(\t)=\Phi_i(\t)^{\ell_i}/\pi_i(\t)^{\ell'_i},\quad 
\gamma_i(\t)=\Phi_i(\t)^{e_i}/\pi_i(\t)^{h_i},
$$
and the identity (\ref{nk}), the second unit may be simplified into:
\begin{align*}
\dfrac{\Phi_i(\t)^{s}\pi_i(\t)^{u+t_k}}{\pi^{n_k}\pi_{i+1}(\t)^t}&\,=\pi^{-n_k}
\Phi_i(\t)^{j-d_ke_i-\ell_it}\pi_i(\t)^{u'_i+jV_i+d_kh_i+t_k+\ell'_it}\\&\,=\pi^{-n_k}
\Phi_i(\t)^{e_i(N+k-d_k)}\pi_i(\t)^{e_0\cdots e_{i-1}n_k-h_i(N+k-d_k)}\\
&\,=\pi^{-n_k}\pi_i(\t)^{e_0\cdots e_{i-1}n_k}\gamma_i(\t)^{N+k-d_k}.
\end{align*}
By (\ref{gamma}), the reduction of the gamma factor is immediate: $$\red_L(\gamma_i(\t)^{N+k-d_k})=z_i^{N+k-d_k}.$$
The unit $\pi^{-1}\pi_i(\t)^{e_0\cdots e_{i-1}}$ depends only on $i$; hence, 
$$\tau_k:=\red_L(\pi^{-n_k}\pi_i(\t)^{e_0\cdots e_{i-1}n_k})=
\red_L(\pi^{-1}\pi_i(\t)^{e_0\cdots e_{i-1}})^{n_k}$$
is a non-zero element in $\F_i$ that depends on $i$, $t$ and $k$.

From (\ref{unit}), (\ref{fraction}) and (\ref{twounits}) we get:
\begin{equation}\label{epsk2}
\epsilon_{k}=\red_L\left(\dfrac{Q_{i,j}(\t)\pi_i(\t)^{t_k}}{\pi^{n_k}\pi_{i+1}(\t)^t}\right)=
R_i(Q_{i,j})(z_i)\cdot\tau_k\cdot z_i^{N+k-d_k}.
\end{equation}
If we use the expression for $R_i(Q_{i,j})$ in Corollary \ref{RQ}, and we consider the element $\zeta_{k}:=(z_{i-1})^{(b_i-j)\ell_{i-1}V_i/e_{i-1}}\cdot\tau_k\in F_i^*$ (that depends on $i$, $t$ and $k$), we get:
\begin{equation}\label{epsk3}
\epsilon_{k}=
\zeta_{k}\cdot z_i^{N}\left(c_{b_i-d_ke_i} z_i^{k-d_k}+c_{b_i-(d_k-1)e_i}z_i^{k-d_k+1}+\cdots+c_{b_i}z_i^k\right)\in\F_{i+1}^*,
\end{equation}

Since $c_{b_i}$ is always non-zero, and $N$ does not depend on $k$, the family of all $\epsilon_{k}$, for $0\le k<f_i$, is an $\F_i$-basis of $\F_{i+1}$. 
\end{proof}

\subsection{A variation on Theorem \ref{Pbasis}}
We can simplify a little bit some quotients and still get an integral basis.
For all $1\le i\le r$, $0\le j<b_i$, define
\begin{equation}\label{Qprime}
Q'_{i,j}(x):=\left\{\begin{array}{ll}
Q_{i,j}(x),&\mbox{ if }j\ne0,\\
1,&\mbox{ if }j=0,
\end{array}
\right.\qquad
H'_{i,j}(x):=\left\{\begin{array}{ll}
H_{i,j}(x),&\mbox{ if }j\ne0,\\
0,&\mbox{ if }j=0,
\end{array}
\right.
\end{equation}
where $Q_{i,j}$, $H_{i,j}$ have still the meaning of Notation \ref{not}.

\begin{thm}\label{Pbasis2}
Let $J=\{(j_0,\dots,j_r)\in\N^{r+1}\mid 0\le j_i<e_if_i,\ \forall\,0\le i\le r\}$, and for each multiindex $\j\in J$ denote:
$$
Q'_\j:=\t^{j_0} Q'_{1,j_1}(\t)\star\cdots\star Q'_{r,j_r}(\t)=\dfrac{\t^{j_0}Q'_{1,j_1}(\t)\cdots Q'_{r,j_r}(\t)}{\pi^{\lfloor H'_{1,j_1}+\cdots +H'_{r,j_r}\rfloor}}\in\B.
$$
Then the family $\bb':=\{Q'_\j\mid \j\in J\}$ is an $\oo$-basis of $\oo_L$ in standard form.
\end{thm}

\begin{proof}
For each $0\le i\le r$, consider the following set: 
$$
B'_i:=\{\t^{j_0} Q'_{1,j_1}\star\cdots\star Q'_{i,j_i}\mid 0\le j_k<e_kf_k,\ \forall\,0\le k\le i\}\subset \B_i.
$$

Denote $\om_j:=Q_{i,j}$, $\om'_j:=Q'_{i,j}$, for all $0\le j<e_if_i$.
Arguing as in the proof of Proposition \ref{PbasisFinal}, the theorem will be proven if we show that the family of all $\om'_j$ satisfies the conditions (a) and (b) of Lemma \ref{valuett}.

Condition (a) is obvious; let us prove condition (b). For fixed $0\le i\le r$, $0\le t<e_0\cdots e_i$, let $0 \le q_t<e_i$ be the integer of Lemma \ref{valuett}, and consider
$$
\as{1.3}
\begin{array}{ll}
 \epsilon_k:=\red_L\left(\pi_{i+1}(\t)^{-t}\om_{q_t+ke_i}\star \pi_i(\t)^{t_k}\right),&\ 0\le k<f_i,\\ \epsilon'_k:=\red_L\left(\pi_{i+1}(\t)^{-t}\om'_{q_t+ke_i}\star \pi_i(\t)^{t'_k}\right), &\ 0\le k<f_i.
\end{array}
$$

If $q_t\ne0$, we have $t'_k=t_k$ and $\epsilon'_k=\epsilon_k$, for all $k$, and we saw along the proof of Proposition \ref{PbasisFinal} that they form an $\F_i$-basis of $\F_{i+1}$.

Suppose $q_t=0$. Then, again, $\epsilon'_k=\epsilon_k$, for all $k\ne0$. For $k=0$ (i.e. $j=0$), let us compute and compare $\epsilon_0$ and $\epsilon'_0$. 
The degree $d_k$ of $R_i(Q_{i,0})(y)$ is zero, because $0\le d_k\le k$; hence, (\ref{epsk3}) shows that: 
$\epsilon_0=\zeta z_i^N$, for some $\zeta\in \F_i^*$.
Since $q_t=0$, $t$ is divisible by $e_i$, and $N=-\ell_it/e_i$; also, since $v(\omega'_j)=0$, we have $t'_0=t/e_i$.
By using again (\ref{ratfracs}) and $\ell_ih_i-\ell'_ie_i=1$, we get:
$$
\pi_{i+1}(\t)^{-e_it'_0}\pi_i(\t)^{t'_0}=\Phi_{i}(\t)^{-\ell_ie_it'_0}\pi_i(\t)^{t'_0+\ell'_ie_it'_0}=\gamma_{i}(\t)^{-\ell_it'_0}=\gamma_{i}(\t)^N.
$$
Hence, $\epsilon'_0=\red_L(\pi_{i+1}(\t)^{-e_it'_0}\pi_i(\t)^{t'_0})=\red_L(\gamma_{i}(\t)^N)=z_i^N$.

Thus, the family $(\epsilon'_k)_{0\le k<f_i}$ differs from the family  $(\epsilon_k)_{0\le k<f_i}$ only in one term: $\epsilon_0=\zeta\epsilon'_0$, for some $\zeta\in \F_i^*$. Since the family $(\epsilon_k)_{0\le k<f_i}$ is an $\F_i$-basis of $\F_{i+1}$, the family $(\epsilon'_k)_{0\le k<f_i}$ has the same property.
\end{proof}

\section{Quotients and global $\mathbf{\p}$-integral bases}\label{secMethod}
We go back to the global context of sections \ref{secOM}, \ref{secMain}, and we keep the notation from those sections.
We denote by $\pset$ the set of prime ideals of $B$ lying above $\p$.  

\subsection{Reduced families of algebraic elements}
For any prime ideal $\P\in\pset$, we use a special notation for two objects attached to the OM representation $\ty_\P$: 
$$
\F_\P:=\F_{r_\P+1,\P}, \qquad \Pi_\P:=\pi_{r_\P+1,\P}(\t)\in L^*.
$$
Thus, $\F_\P$ is a computational representation of the local residue field $\oo_\P/\P\oo_\P$. The rational fractions $\pi_{i,\P}$ were defined in (\ref{ratfracs}); recall that $w_\P(\Pi_\P)=1/e(\P/\p)$.

The concept of a \emph{$\p$-reduced fa\-mily} of elements of the function field of a curve was introduced by W. M. Schmidt in the context of Puiseux expansions \cite{schmidt,schoernig}. In this section we use these ideas conveniently adapted to our more general setting.

\begin{definition}
Consider the following $\p$-valuation mapping:
$$
w:=w_\p\colon L\lra \Q\cup\{\infty\},\quad w(\alpha)=\min\nolimits_{\P\in\pset}\{w_\P(\alpha)\}.
$$
\end{definition}

Clearly, $w(a)=v_\p(a)$, for all $a\in K$. The map $w$ does not behave well with respect to multiplication, but it  has some of the typical properties of a valuation.

\begin{lemma}\label{valuation2}
Let $a\in K$, and $\alpha,\beta\in L$.
\begin {enumerate}
\item $w(a\alpha)=w(a)+w(\alpha)=v_\p(a)+w(\alpha)$.
\item $w(\alpha+\beta)\ge \min\{w(\alpha),w(\beta)\}$, and if $w(\alpha)\ne w(\beta)$, then equality holds.
\end {enumerate}
\end{lemma}

\begin{definition}
For any value $\delta\in w(L)$, we denote $$L_\delta:=\{\alpha\in L\mid w(\alpha)\ge \delta\}\supset L_\delta^+:=\{\alpha\in L\mid w(\alpha)>\delta\}.$$ 
\end{definition}

Note that $L_\delta\subset B_\p$ if $\delta\ge0$. These subgroups $L_\delta$, $L_\delta^+$ have a natural structure of $A_\p$-modules. Since $\p L_\delta\subset L_\delta^+$, the quotient $L_\delta/L_\delta^+$ has a natural structure of $\F_\p$-vector space.

\begin{definition}
Consider the $\F_\p$-vector space, $V:=\prod_{\P\in\pset}\F_\P$, of dimension $\sum_{\P\in\pset}f(\P/\p)$.
For each $\delta\in w(L)$, $\delta\ge 0$, we define a kind of reduction map:
$$
\red_\delta\colon L_\delta\lra V,\quad \red_\delta(\alpha)=(\alpha_{\delta,\P})_{\P\in\pset},\quad
\alpha_{\delta,\P}=
\red_{L_\P}\left( i_\P\left(\alpha/\Pi_\P^{e(\P/\p)\delta}\right)\right).
$$
\end{definition}

Note that $\alpha_{\delta,\P}=0$ if and only if $w_\P(\alpha)>\delta$. 
Clearly, $\red_\delta$ is an homomorphism of $A_\p$-modules, and $\ker(\red_\delta)=L_\delta^+$. Therefore, $\red_\delta$ induces an embedding of $L_\delta/L_\delta^+$ as an $\F_\p$-subspace of $V$.

\begin{definition}
A finite subset $\bb=\{\alpha_1,\dots,\alpha_m\}\subset L$ is called \emph{$\p$-reduced} if for all families $a_1,\dots,a_m\in A_\p$, one has:
\begin{equation}\label{reduceness}
w\left(\sum\nolimits_{1\le i\le m}a_i\alpha_i\right)=\min\{w(a_i\alpha_i)\mid 1\le i\le m\}.
\end{equation}
\end{definition}

Let $\nu:=\min_{1\le i \le m}\{v_\p(a_i)\}$.
The left and right terms of (\ref{reduceness}) diminish both by $\nu$ if we replace all $a_i$ by $a'_i=a_i/\pi^\nu$. Thus, in order to check the equality (\ref{reduceness}) we can always assume that not all elements $a_1,\dots, a_m\in A_\p$ belong to $\p A_\p$.   

Reduceness is a sufficient condition to ensure that certain families of elements in $B_\p$ are a $\p$-integral basis.

\begin{lemma}\label{reducedbasis}
For $n=[L\colon K]$, let $\bb=\{\alpha_1,\dots,\alpha_n\}\subset L$ be a $\p$-reduced set such that $0\le w(\alpha)<1$, for all $\alpha\in\bb$. Then, $\bb$ is a $\p$-integral basis of $B/A$.
\end{lemma}

\begin{proof}
We need only to check that the elements of $\bb$ are linearly independent modulo $\p B_\p$. Suppose $\sum_{1\le i\le n}a_i\alpha_i\in \p B_\p$, for certain $a_1,\dots, a_n\in A_\p$. Since $w\left(\sum_{1\le i\le n}a_i\alpha_i\right)\ge 1$ and $\bb$ is reduced, we have $w(a_i\alpha_i)\ge 1$, for all $1\le i\le n$. Since $w(\alpha_i)<1$, this implies that $w(a_i)>0$, or equivalently, $a_i\in \p A_\p$, for all $i$.  
\end{proof}

Let us give a more practical criterion to check that a subset of $L$ is reduced.

\begin{lemma}\label{criterion}
Let $\bb\subset L$ be a finite subset such that $0\le w(\alpha)<1$, for all $\alpha\in\bb$. For each $\delta\in w(L)$, denote $\bb_\delta:=\{\alpha\in \bb\mid w(\alpha)=\delta\}$. Then, $\bb$ is $\p$-reduced if and only if $\red_\delta(\bb_\delta)$ is an $\F_\p$-linearly independent family of $V$, for all $\delta\in w(\bb)$.
\end{lemma}

\begin{proof}
Write $\bb=\{\alpha_1,\dots,\alpha_m\}$, and let $I_\delta:=\{1\le i\le m\mid \alpha_i\in\bb_\delta\}$, for each $\delta\in w(\bb)$. For any family $(a_i)_{i\in I_\delta}$ of elements in $A_\p$, we clearly have:
\begin{equation}\label{conversion}
\red_\delta\left(\sum\nolimits_{i\in I_\delta}a_i\alpha_i\right)=
\sum\nolimits_{i\in I_\delta}\red_\delta(a_i\alpha_i)=
\sum\nolimits_{i\in I_\delta}\overline{a_i}\red_\delta(\alpha_i),
\end{equation}
where $\overline{a_i}\in\F_\p$ is the class of $a_i$ modulo $\p A_\p$.

Suppose $\bb$ is a reduced set. If $\overline{a_{i_0}}\ne0$, for some $i_0\in I_\delta$, then $w(a_{i_0}\alpha_{i_0})=\delta=\min_{i\in I_\delta}\{w(a_i\alpha_i)\}$. By reduceness, we get $w\left(\sum_{i\in I_\delta}a_i\alpha_i\right)=\delta$, and (\ref{conversion}) shows that $\sum_{i\in I_\delta}\overline{a_i}\red_\delta(\alpha_i)\ne0$. Thus, the fa\-mily $\red_\delta(\bb_\delta)$ is $\F_\p$-linearly independent.

Conversely, suppose that $\red_\delta(\bb_\delta)$ is an $\F_\p$-linearly independent family of $V$, for all $\delta\in w(\bb)$. Take $a_1,\dots,a_m\in A_\p$, not all of them belonging to $\p A_\p$, and let
$$\delta:=\min\{w(a_i\alpha_i)\mid 1\le i\le m\},\quad J_\delta:=\{1\le i\le m\mid w(a_i\alpha_i)=\delta\}.
$$
Since $0\le w(\alpha)<1$, for all $\alpha\in\bb$, and not all $a_m$ belong to $\p A_\p$, we have $\delta<1$. Thus, $J_\delta\subset I_\delta$, and (\ref{conversion}) shows that $w\left(\sum_{i\in J_\delta}a_i\alpha_i\right)=\delta$. Since $w\left(\sum_{i\not\in J_\delta}a_i\alpha_i\right)>\delta$, we get $w\left(\sum_{i=1}^ma_i\alpha_i\right)=\delta$, as desired.  
\end{proof}

\subsection{Domination and similarity of prime ideals}\label{secDomination}
The aim of the Montes algorithm is to determine successive dissections of the set $\pset$, till each prime ideal lying over $\p$ is singled out. In this section, we derive from these dissections a partial ordering on a quotient of $\pset$ by a certain equivalence relation. Needless to say, these relationships between prime ideals are not intrinsic; they depend on the choice of the polynomial $f(x)\in A[x]$, defining the extension $L/K$. For instance, the first dissection of $\pset$ is determined by the factorization of $f(x)$ modulo $\p$:
$$
\overline{f}(y)=\prod\nolimits_{\varphi}\varphi(y)^{a_{\varphi}},
$$
into a product of powers of pairwise different monic irreducible polynomials $\varphi\in\F_\p[y]$. By Hensel's lemma, this determines a partition 
$$
\pset=\coprod\nolimits_\varphi \pset_\varphi,\quad \pset_\varphi:=\left\{\P\in\pset\mid \overline{F}_\P\mbox{ is a power of }\varphi\right\}.
$$

Let us briefly recall how the Montes algorithm proceeds to obtain further dissections of the subsets $\pset_\varphi$. We use the version of the algorithm described in \cite[Sec. 4]{BNS}, guaranteeing that the OM representations $\ty_\P$ have order $r_\P+1$, where $r_\P$ is the Okutsu depth of $F_\P$.

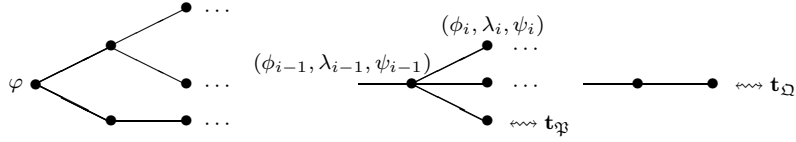
\begin{figure}\caption{Tree $\ttt_\varphi$ of OM representations of the irreducible factors of $f(x)$ whose reduction modulo $\p$ is a power of $\varphi$.}
\label{figTree}
\begin{center}
\setlength{\unitlength}{5.mm}
\begin{picture}(20,3.4)
\put(-.15,.6){$\bullet$}\put(1.85,1.65){$\bullet$}
\put(1.85,-.35){$\bullet$}
\put(0,.8){\line(2,1){2}}\put(0.02,.8){\line(2,1){2}}
\put(0,.8){\line(2,-1){2}}\put(0.02,.8){\line(2,-1){2}}
\put(-.7,.7){\begin{footnotesize}$\varphi$\end{footnotesize}}
\put(3.85,2.65){$\bullet$}
\put(3.85,.65){$\bullet$}
\put(2,1.8){\line(2,1){2}}\put(0.02,.8){\line(2,1){2}}
\put(2,1.8){\line(2,-1){2}}\put(0.02,.8){\line(2,-1){2}}
\put(2,-.2){\line(2,0){2}}\put(0.02,.8){\line(2,-1){2}}
\put(3.85,-.35){$\bullet$}
\put(4.5,.6){\begin{footnotesize}$\cdots$\end{footnotesize}}
\put(4.5,-.4){\begin{footnotesize}$\cdots$\end{footnotesize}}
\put(4.5,2.6){\begin{footnotesize}$\cdots$\end{footnotesize}}
\put(9.85,.65){$\bullet$}\put(11.85,1.65){$\bullet$}
\put(11.85,-.35){$\bullet$}\put(11.85,.65){$\bullet$}
\put(8.6,.8){\line(1,0){3.4}}
\put(10,.8){\line(2,1){2}}\put(10.02,.8){\line(2,1){2}}
\put(10,.8){\line(2,-1){2}}\put(10.02,.8){\line(2,-1){2}}
\put(5.8,1.2){\begin{footnotesize}$(\phi_{i-1},\lambda_{i-1},\psi_{i-1})$\end{footnotesize}}
\put(10.8,2.2){\begin{footnotesize}$(\phi_i,\lambda_i,\psi_i)$\end{footnotesize}}
\put(12.6,-.4){\begin{footnotesize}$ \leftrightsquigarrow\ty_{\P} $\end{footnotesize}}
\put(12.7,.6){\begin{footnotesize}$\cdots$\end{footnotesize}}
\put(12.7,1.6){\begin{footnotesize}$\cdots$\end{footnotesize}}
\put(18.6,.6){\begin{footnotesize}$\leftrightsquigarrow\ty_{\q} $\end{footnotesize}}
\put(14.6,.8){\line(1,0){3.4}}
\put(15.85,.65){$\bullet$}\put(17.85,.65){$\bullet$}
\end{picture}
\end{center}
\end{figure}

For each $\varphi$, we consider initially a triple $(\ty,\phi,\omega)$, where  $\ty=(\varphi)$ is a type of order zero, $\phi$ is a representative of $\ty$ (a monic lift of $\varphi$ to $A[x]$) and $\omega=\ord_\varphi(\overline{f})$. We submit this triple to a kind of \emph{branching process}, by enlarging $\ty$ to different types of higher order. This process is repeated for each branch, till all OM representations of the prime ideals in $\pset_\varphi$ are obtained. We build in this way a connected tree $\ttt_\varphi$ of OM representations, whose root node is labelled by the polynomial $\varphi$, and the rest of the nodes are labelled by triples $(\phi,\lambda,\psi)$. The prime ideals of $\pset_\varphi$ are in 1-1 correspondence with the leaves of the tree, and the type $\ty_\P$ attached to a leaf is obtained by gathering the invariants of all nodes in the unique path joining the leaf to its root node (see Figure \ref{figTree}).

In a general iteration, the branching process is applied to a triple $(\ty,\phi,\omega)$, where $\ty$ is a strongly optimal type of order $i-1\ge 0$, dividing $f(x)$, $\phi$ is a representative of $\ty$ and $\omega$ a positive integer. We compute the Newton polygon $N_{i,\omega}(f)\subset N_{\phi,v_i}^ -(f)$ determined by the first $\omega+1$ coefficients of the $\phi$-expansion of $f(x)$. The branches of $\ty$ are determined by all pairs $(\lambda,\psi)$, where $\lambda$ runs on all slopes of the sides of $N_{i,\omega}(f)$, and for each $\lambda$, the polynomial $\psi$ runs on all monic irreducible factors of $R_{\lambda,i}(f)$. 

If $\omega=1$, there is only one branch, and the triple $(\phi,\lambda,\psi)$ determines a leaf of $\ttt_\varphi$. If  $\omega>1$, for each branch $(\lambda,\psi)$ we consider the type  $\ty_{\lambda,\psi}:=(\ty;(\phi,\lambda,\psi))$, of order $i$, we compute a representative $\phi_{\lambda,\psi}$ of this type, and the positive integer $\omega_{\lambda,\psi}:=\ord_{\ty_{\lambda,\psi}}(f)$. The following subsets of $\pset_\varphi$:
$$
\pset_{\lambda,\psi}:=\left\{\P\in\pset_\varphi \mid  \ty_{\lambda,\psi}\mbox{ divides }F_\P\right\}
$$
are pairwise disjoint. If $\ty_{\lambda,\psi}$ is strongly optimal, then $(\phi,\lambda,\psi)$ labels a new node of $\ttt_\varphi$, and we submit the triple $(\ty_{\lambda,\psi},\phi_{\lambda,\psi},\omega_{\lambda,\psi})$ to further branching. Otherwise, $\phi_{\lambda,\psi}$ is also a representative of $\ty$, and we submit the triple  $(\ty,\phi_{\lambda,\psi},\omega_{\lambda,\psi})$ to further branching; this is called a \emph{refinement step}. 

\begin{definition}
Let $\pset_0:=\bigcup_{\ord_\varphi(\overline{f})=1}\pset_\varphi$ be the subset of $\pset$ formed by the prime ideals singled out by the first dissection.
\end{definition}

In fact, if $\ord_\varphi(\overline{f})=1$, then $\pset_\varphi=\{\P\}$ consists of a single prime ideal, with $e(\P/\p)=1$, $f(\P/\p)=\deg \varphi$. In the initial step, we have already $\omega=1$, so that $F_\P$ has Okutsu depth zero and the tree $\ttt_\varphi$ has only the root node and one leaf:
\begin{center}
\setlength{\unitlength}{5.mm}
\begin{picture}(10,1)
\put(-.15,0){$\bullet$}\put(1.85,0){$\bullet$}
\put(0,0.15){\line(1,0){2}}\put(0,0.17){\line(1,0){2}}
\put(-.7,0.1){\begin{footnotesize}$\varphi$\end{footnotesize}}
\put(2.5,0.1){\begin{footnotesize}$(\phi,\lambda,\psi)$\end{footnotesize}}
\put(4.4,0.1){\begin{footnotesize}$\quad \leftrightsquigarrow\ \ty_{\P}=(\varphi;(\phi,\lambda,\psi))$\end{footnotesize}}
\end{picture}
\end{center}

\begin{definition}
A leaf of $\ttt_\varphi$ is \emph{isolated} if it is the unique branch of its previous node. 
We say that $\P\in\pset$ is \emph{isolated} if the leaf corresponding to $\ty_\P$ is isolated.
\end{definition}

For instance, in Figure \ref{figTree}, the prime $\P$ is non-isolated and the prime $\q$ is isolated.

\begin{definition}
Let $\P,\q\in\pset$, $\P\ne\q$. If $\psi_{0,\P}=\psi_{0,\q}$ (that is, $\ty_\P$ and $\ty_\q$ belong to the same connected tree of OM representations), we define the \emph{index of coincidence} between $\ty_\P$
and $\ty_\q$ as:
$$
i(\ty_\P,\ty_\q)=
\min\left\{j\in\Z_{>0}\tq (\phi_{j,\P},\lambda_{j,\P},\psi_{j,\P})\ne
(\phi_{j,\q},\lambda_{j,\q},\psi_{j,\q})\right\}.
$$
If $\psi_{0,\P}\ne\psi_{0,\q}$, we define $i(\ty_\P,\ty_\q)=0$.  
\end{definition}

By the very definition, we have:
\begin{equation}\label{tope}
i(\ty_\P,\ty_\q)\le
\begin{cases}
r_\P,&\mbox{if $\P$ is isolated},\\r_\P+1,&\mbox{if $\P$ is non-isolated.}
\end{cases}
\end{equation}

If $j=i(\ty_\P,\ty_\q)$, the types $\ty_\P$, $\ty_\q$ have the same truncation at the $(j-1)$-th order:
$\ty_{j-1}:=\op{Trunc}_{j-1}(\ty_\P)=\op{Trunc}_{j-1}(\ty_\q)$.
The last level of this type labels the first node of $\ttt_\varphi$, where the branches corresponding to the leaves $\ty_\P$, $\ty_\q$ diverge. 
 
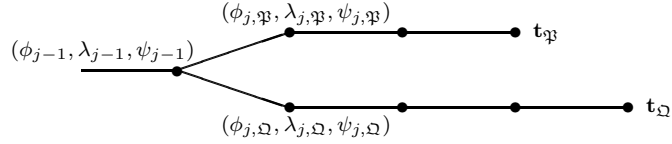
\begin{center}
\setlength{\unitlength}{5.mm}
\begin{picture}(12,4)
\put(.85,1.6){$\bullet$}\put(3.85,2.65){$\bullet$}\put(3.85,.65){$\bullet$}
\put(-1.5,1.8){\line(1,0){2.5}}
\put(1,1.8){\line(3,1){3}}\put(1,1.83){\line(3,1){3}}
\put(1,1.8){\line(3,-1){3}}\put(1,1.83){\line(3,-1){3}}
\put(6.85,2.65){$\bullet$}\put(6.85,.65){$\bullet$}
\put(4,2.8){\line(1,0){3}}\put(4,2.83){\line(1,0){3}}
\put(4,.8){\line(1,0){3}}\put(4,.83){\line(1,0){3}}
\put(9.85,2.65){$\bullet$}\put(9.85,.65){$\bullet$}
\put(7,2.8){\line(1,0){3}}\put(7,2.83){\line(1,0){3}}
\put(7,0.8){\line(1,0){3}}\put(7,.83){\line(1,0){3}}
\put(12.85,.65){$\bullet$}
\put(10,0.8){\line(1,0){3}}\put(10,.83){\line(1,0){3}}
\put(-3.4,2.2){\begin{footnotesize}$(\phi_{j-1},\lambda_{j-1},\psi_{j-1})$\end{footnotesize}}
\put(10.5,2.7){\begin{footnotesize}$\ty_\P$\end{footnotesize}}
\put(13.5,.7){\begin{footnotesize}$\ty_\q$\end{footnotesize}}
\put(2.2,3.2){\begin{footnotesize}$(\phi_{j,\P},\lambda_{j,\P},\psi_{j,\P})$\end{footnotesize}}
\put(2.2,.2){\begin{footnotesize}$(\phi_{j,\q},\lambda_{j,\q},\psi_{j,\q})$\end{footnotesize}}
\end{picture}
\end{center}

The polynomials $\phi_{j,\P}$, $\phi_{j,\q}$ are representatives of $\ty_{j-1}$, but they do not necessari\-ly coincide. 
Nevertheless, there exists a \emph{greatest common $\phi$-polynomial} $\phi(\P,\q)$ of the pair $\ty_\P,\,\ty_\q$ \cite[Defn. 3.7]{newapp}. The algorithm computes at some iteration a representative $\phi(\P,\q)$ of $\ty_{j-1}$, admitting two different branches, $(\lambda_\P^\q,\psi_\P^\q)\ne(\lambda_\q^\P,\psi_\q^\P)$, leading, eventually after some refinement steps and/or further branching, to the nodes  $(\phi_{j,\P},\lambda_{j,\P},\psi_{j,\P})$, $(\phi_{j,\q},\lambda_{j,\q},\psi_{j,\q})$, respectively.   
The slopes $\lambda_\P^\q,\,\lambda_\q^\P$ are called the \emph{hidden slopes} of the pair $\ty_\P,\,\ty_\q$. 

Let $(\ty,\phi_r,\omega)$ be one of the triples submitted to the branching process, along the flow of the Montes algorithm, and suppose that $\omega>1$. Recall that $\ty$ is a strongly optimal type of order (say) $r-1\ge0$, and $\phi_r$ is a representative of $\ty$. Let $S$ be a side of $N_{r,\omega}(f)\subset N_r^-(f)$, $\lambda\in\Q^-$ the slope of $S$, and 
$$
R_{\lambda,r}(f)(y)\sim \psi_{1}(y)^{n_1}\cdots\psi_{t}(y)^{n_t}
$$
the factorization of $R_{\lambda,r}(f)(y)$ into the product of pairwise different monic irreducible polynomials in $\F_r[y]$. 
Write $\lambda=-h/e$, with $h,e$ positive coprime integers. The length $\ell(S)$ of the side $S$ is, by definition,  the length of the projection of $S$ to the horizontal axis. By Definition \ref{respol},  
$$\ell(S)=e\deg R_{\lambda,r}(f)=e\sum\nolimits_{1\le k\le t}n_k\deg \psi_k.
$$ 

\begin{definition}
We define the \emph{terminal length} of $S$ as:
$$
\lt(S):=e\sum\nolimits_{n_k=1}\deg\psi_k.
$$ 
We say that $S$ is a \emph{terminal side of order $r$}, if $\lt(S)>0$, or equivalently, if at least one irreducible factor of $R_{\lambda,r}(f)$ divides this polynomial with exponent one.
\end{definition}

Let $S$ be a terminal side of order $r$. Each branch $(\lambda,\psi)$ of $(\ty,\phi,\omega)$, with 
$$
\omega_{\lambda,\psi}:=\ord_{\ty_{\lambda,\psi}}(f)=\ord_{\psi}(R_{\lambda,r}(f))=1,
$$ singles out a prime ideal $\P_{\lambda,\psi}$ of $\pset$. In fact, 
$1=\ord_{\ty_{\lambda,\psi}}(f)=\sum\nolimits_{\P\in\pset}\ord_{\ty_{\lambda,\psi}}(F_\P)$,
so that $\ord_{\ty_{\lambda,\psi}}(F_\P)=0$, for all $\P\in\pset$, except for one prime ideal, say $\P_{\lambda,\psi}$, for which $\ord_{\ty_{\lambda,\psi}}(F_{\P_{\lambda,\psi}})=1$. We denote 
$$
\pset_S:=\{\P_{\lambda,\psi}\mid\ord_{\psi}(R_{\lambda,r}(f))=1\}\subset\pset_\varphi.
$$
Note that $\ty_{\lambda,\psi}$ is simultaneously $f$-complete and $F_{\P_{\lambda,\psi}}$-complete.

The representative $\phi_{\lambda,\psi}$ of $\ty_{\lambda,\psi}$ has degree $e(\deg\psi) m_r$. Hence, $\ty_{\lambda,\psi}$ is strongly optimal if and only if $e\deg\psi>1$; in this case, $F_{\P_{\lambda,\psi}}$ has Okutsu depth $r_{\P_{\lambda,\psi}}=r$, and $(\phi_r,\lambda,\psi)$ is the $r$-th level of the OM representation of $\P_{\lambda,\psi}$. The prime ideal $\P_{\lambda,\psi}$ is isolated, because $\ty_{\lambda,\psi}$ has a unique branch, which is a leave of the tree $\ttt_\varphi$. This $(r+1)$-th level of $\ty_{\P_{\lambda,\psi}}$ is constructed by a last iteration applied to the triple $(\ty_{\lambda,\psi},\phi_{\lambda,\psi},\omega_{\lambda,\psi}=1)$.  

If $e\deg\psi=1$, the Okutsu depth of $F_{\P_{\lambda,\psi}}$ is $r_{\P_{\lambda,\psi}}=r-1$, and $\P$ is non-isolated. In fact, the iteration that constructs the leaf attached to $\P_{\lambda,\psi}$ is applied to the triple $(\ty,\phi_{\lambda,\psi},\omega_{\lambda,\psi}=1)$; hence, it yields a node of level $r$ of the tree of OM representations. On the other hand, since $\omega>1$, the initial triple $(\ty,\phi_r,\omega)$ has other branches of level $r$.

\begin{definition}\label{partition}
Let $\st$ be the set of all terminal sides that occur along the application of the Montes algorithm to $f(x)$ and $\p$. 
We have a partition: 
$$
\pset=\pset_0\cup\bigcup\nolimits_{S\in\st}\pset_S.
$$

Let $S\in\st$ be a terminal side of order $r$. If $\P\in\pset_S$, we denote $\phi_S:=\phi_r$, the $\phi$-polynomial from which the side $S$ was derived. Note that $\phi_S=\phi_{r,\P}$, if $\P$ is isolated, but $\phi_S$ is not a $\phi$-polynomial of $\ty_\P$, if $\P$ is non-isolated.     
\end{definition}

\begin{definition}\label{domination}
Let $S\in\st$, and let $\P\in\pset_S$. We say that $\q\in\pset\setminus\pset_0$ \emph{dominates} $\P$, and we write $\q\succ\P$, if $w_\q(\phi_S(\t))\ge w_\P(\phi_S(\t))$.

We say that $\P,\q$ are \emph{similar}, and we write $\P\simeq \q$,  if $\q\succ\P$ and $\P\succ\q$.
\end{definition}

\begin{lemma}\label{criteria}
Let $S\in\st$ have order $r$. Let $\P\in\pset_S$ and $\q\in\pset\setminus\pset_0$, $\q\ne\P$.
\begin{enumerate}
\item $\q\succ\P$ if and only if $i(\ty_\P,\ty_\q)=r$, $\phi_S=\phi(\P,\q)$ and $\,|\lambda_\q^\P|\ge |\lambda_\P^\q|$.
\item $\q\simeq \P$ if and only if $\q\in\pset_S$.
\end{enumerate}
\end{lemma}

\begin{proof}
Let $\lambda$ be the slope of $S$. By \cite[Thm. 3.1]{HN}, 
\begin{equation}\label{knownP}
w_\P(\phi_S(\t))=\left(V_{r,\P}+|\lambda|\right)/e_{0,\P}\cdots e_{r-1,\P}.
\end{equation}

As mentioned above, $r=r_\P$, if $\P$ is isolated, and $r=r_\P+1$, otherwise. By (\ref{tope}), $j:=i(\ty_\P,\ty_\q)\le r$. 
If $j=0$, we have $w_\q(\phi_S(\t))=0$, and the claimed equivalent conditions of item 1 are both false. If $j>0$,  \cite[Prop. 3.8]{newapp} shows that
\begin{equation}\label{knownQ}
w_\q(\phi_S(\t))=\begin{cases}
\dfrac{V_{r,\P}+|\lambda_\q^\P|}{e_{0,\P}\cdots e_{r-1,\P}},&\mbox{ if $j=r$ and }\phi(\P,\q)=\phi_S,\\
\dfrac{m_{r,\P}}{m_j}\,\dfrac{V_{j}+\min\{|\lambda_\P^\q|,|\lambda_\q^\P|\}}{e_0\cdots e_{j-1}},&\mbox{ otherwise}.
\end{cases}
\end{equation}
The Okutsu invariants $e_0,\dots,e_{j-1},m_j,V_j$ of $\ty_\P$ and $\ty_\q$ coincide, and for them we dropped the subindex $\P$, or $\q$. 

Suppose $j=r$ and $\phi(\P,\q))=\phi_S$. We then have $\lambda_\P^\q=\lambda$, by the definition of the hidden slope. By (\ref{knownP}) and (\ref{knownQ}),
$\,|\lambda_\q^\P|\ge |\lambda_\P^\q|$ is equivalent to $\q\succ\P$. 
Therefore, in order to prove item 1 of the lemma, it is sufficient to check that $\q\succ\P$ implies  $j=r$ and $\phi(\P,\q))=\phi_S$.

In every refinement step, the slope grows strictly in absolute size \cite[Thm. 3.1]{bordeaux}. Hence,
\begin{equation}\label{refinement}
|\lambda_\P^\q|\le |\lambda_{j,\P}|,\ \mbox{ and }\ \phi(\P,\q)\ne\phi_S\,\Longrightarrow\, |\lambda_\P^\q|<|\lambda|.
\end{equation}

If $j=r$ and $\phi(\P,\q)\ne\phi_S$, then we get directly $w_\q(\phi_S(\t))<w_\P(\phi_S(\t))$, by (\ref{knownP}), (\ref{knownQ}) and (\ref{refinement}). 
If $j<r$, then (\ref{knownP}), (\ref{knownQ}), (\ref{refinement}) and the explicit recurrent formulas for $V_i$ from section \ref{secOkutsu}, show that
\begin{align*}
w_\q(\phi_S(\t))\le&\ m_{r,\P}\,\dfrac{V_{j}+|\lambda_{j,\P}|}{m_j\,e_0\cdots e_{j-1}}
= m_{r,\P}\,\dfrac{V_{j+1,\P}}{m_{j+1,\P}e_{0,\P}\cdots e_{j,\P}}\\\le &\ 
m_{r,\P}\,\dfrac{V_{r,\P}}{m_{r,\P}\,e_{0,\P}\cdots e_{r-1,\P}}<w_\P(\phi_S(\t)).
\end{align*}
This ends the proof of item 1.

Let us now prove item 2. If $\q\in\pset_S$, we have by construction: $i(\ty_\P,\ty_\q)=r$, $\phi(\P,\q)=\phi_S$ and $\lambda_\P^\q=\lambda_\q^\P=\lambda$. Hence, $\q\succ\P$ and $\P\succ\q$, by the first item. 

Conversely, suppose $\P\simeq\q$. Let $T$ be the terminal side for which  $\q\in\pset_T$, and let $\mu$ be the slope of $T$. By the first item, $i(\ty_\P,\ty_\q)=r$, $\phi(\P,\q)=\phi_S=\phi_T$, and $\lambda=\lambda_\P^\q= \lambda_\q^\P=\mu$. Hence, $S=T$.
\end{proof}

\begin{lemma}\label{partorder}
The relation of domination is reflexive and transitive. Thus, it induces a partial ordering on the set $(\pset\setminus\pset_0)/\!\simeq\,$ of  similarity classes of $\pset\setminus\pset_0$. 
\end{lemma}

\begin{proof}
The reflexive property is obvious. Let us prove transitivity. Suppose $\l\succ\q$, $\q\succ\P$. Let $S,T\in\st$ be the terminal sides such that $\P\in\pset_S$, $\q\in\pset_T$. By Lemma \ref{criteria}, $\phi_S=\phi(\P,\q)$, $\phi_T=\phi(\q,\l)$, $r:=i(\ty_\P,\ty_\q)$ is equal to $r_\P$, or $r_\P+1$, according to $\P$ being isolated or not, and $s:=i(\ty_\q,\ty_\l)$ is equal to $r_\q$, or $r_\q+1$, according to $\q$ being isolated or not. By (\ref{tope}), we have $r\le s$, so that $i(\ty_\P,\ty_\l)=r$. 

Let $(\ty,\phi,\omega)$ be the first triple such that the three prime ideals $\P,\q,\l$ do not belong to the same of its branches. Let $(\lambda,\psi)$ be the branch to which $\P$ belongs; that is, $\P\in\pset_{\lambda,\psi}$. The prime ideal $\q$ cannot belong to the same branch. In fact, this would separate $\q$ from $\l$, and we would have $\phi=\phi(\q,\l)=\phi_T$; but this is impossible, because the branch of $\phi_T$ to which $\q$ belongs contains no other prime ideal. Therefore, $\phi=\phi(\P,\q)=\phi_S$; in particular, $\pset_{\lambda,\psi}=\{\P\}$. Hence, $\P$ and $\l$ are also separated by this triple, and this implies $\phi(\P,\l)=\phi=\phi_S$.

By Lemma \ref{criteria}, in order to prove that $\l\succ\P$, we need only to show that $|\lambda_\l^\P|\ge |\lambda_\P^\l|$. We now have two possibilities, according to $\l$, $\q$ belonging to the same branch, or to different branches. 

\begin{center}
\setlength{\unitlength}{4.mm}
\begin{picture}(22,3)
\put(0.2,1.3){\line(1,1){1}}\put(0.2,1.3){\line(1,-1){1}}
\put(-2.5,1.2){\begin{footnotesize}$\phi(\P,\q)$\end{footnotesize}}
\put(1.4,2.2){\begin{footnotesize}$(\lambda,\psi)\ \ \quad \pset_{\lambda,\psi}=\{\P\}$\end{footnotesize}}
\put(1.3,0){\begin{footnotesize}$(\lambda',\psi')\ \quad \q,\,\l\in\pset_{\lambda',\psi'}  $\end{footnotesize}}
\put(12.5,1.2){\begin{footnotesize}$\phi(\P,\q)$\end{footnotesize}}
\put(15.2,1.3){\line(1,1){1}}\put(15.2,1.3){\line(1,-1){1}}\put(15.2,1.3){\line(1,0){1}}
\put(16.6,2.3){\begin{footnotesize}$(\lambda,\psi)\quad\,  \quad \pset_{\lambda,\psi}=\{\P\}$\end{footnotesize}}
\put(16.5,1.2){\begin{footnotesize}$(\lambda',\psi')\ \ \quad \pset_{\lambda',\psi'}=\{\q\}$\end{footnotesize}}
\put(16.4,0){\begin{footnotesize}$(\lambda'',\psi'')\ \quad \l\in\pset_{\lambda'',\psi''}$\end{footnotesize}}
\end{picture}
\end{center}

In the first case, we have $|\lambda_\l^\P|=|\lambda'|=|\lambda_\q^\P|\ge|\lambda_\P^\q|=|\lambda|=|\lambda_\P^\l|$. In the second case, the argument is similar: $\,|\lambda_\l^\P|=|\lambda_\l^\q|\ge|\lambda_\q^\l|= |\lambda_\q^\P|\ge |\lambda_\P^\q|=|\lambda_\P^\l|$.
\end{proof}

By Lemma \ref{criteria}, there is a natural bijection between $\st$ and $(\pset\setminus\pset_0)/\simeq$. Therefore, domination induces a partial ordering on $\st$ as well.

\subsection{Method of the quotients}
For each $\P\in\pset_0$, with OM representation $\ty_\P=(\psi_{0,\P};(\phi_{1,\P},\lambda_{1,\P},\psi_{1,\P}))$, denote by $Q_\P(x)$ the quotient of the division with remainder of $f(x)$ by $\phi_{1,\P}(x)$. Consider the set:
$$
\bb_{\pset_0}:=\bigcup\nolimits_{\P\in\pset_0}\bb_\P,\quad \bb_\P:=\{Q_\P(\t),\t Q_\P(\t),\dots,\t^{f_{0,\P}-1}Q_\P(\t)\}. 
$$

Let $S$ be a terminal side of order $r$, derived from a type 
\begin{equation}\label{typeS}
\ty=(\psi_0;(\phi_{1},\lambda_{1},\psi_{1});\cdots;(\phi_{r-1},\lambda_{r-1},\psi_{r-1})),
\end{equation}
with representative $\phi_r$. Denote by $\lambda_r$ the slope of $S$. For all $1\le i\le r$, let $b_i$ be the abscissa of the right end point of the side of slope $\la_i$ of $N_i^-(f)$. 
For all $0\le j<b_i$, let $Q_{i,j}$ be the $(b_i-j)$-th quotient of the $\phi_i$-expansion of $f(x)$. Consider the set:
$$
J_S:=\{(j_0,\dots,j_{r-1},j)\in\N^{r+1}\mid 0\le j_i<e_if_i, \ 0\le i<r; \ 0\le j<\lt(S)\},
$$ 
and for any $\j\in J_S$, consider the element: 
\begin{equation}\label{generic}
Q_\j:=\dfrac{\t^{j_0} Q'_{1,j_1}(\t)\cdots Q'_{r-1,j_{r-1}}(\t)\, Q_{r,j}(\t)}{\pi^{\lfloor H'_{1,j_1}+\cdots +H'_{r-1,j_{r-1}}+H_{r,j} \rfloor}}\in B_\p,
\end{equation}
where $Q'_{i,j}$, $H'_{i,j}$, $H_{i,j}$, are defined in (\ref{Qprime}). Finally, let
 $\bb_S:=\{Q_\j\mid \j\in J_S\}$.

\begin{thm}\label{pBasis2}
The following family is a $\p$-reduced $\p$-integral basis of $B/A$:
$$
\bb:=\bb_{\pset_0}\cup\left(\bigcup\nolimits_{S\in \st} \bb_S\right).
$$
\end{thm}

For the proof of the theorem we need two lemmas.

\begin{lemma}\label{auxiliar}
Let $S$ be a terminal side of order $r$, derived from a type $\ty$ with representative $\phi_r$, as in (\ref{typeS}). Let $\lambda_r$ be the slope of $S$. For each $\P\in\pset_S$, denote by $\psi_\P\in\F_r[y]$ the irreducible factor of $R_r(f)$, such that $\P=\P_{\lambda_r,\psi_\P}$ is determined by the branch $(\lambda_r,\psi_\P)$.
 
\begin{enumerate}
\item For each $0\le j<b_r$, $w_\P(Q_{r,j}(\t))=H_{r,j}$ if and only if $\psi_\P\nmid R_r(Q_{r,j})$. 

If $0\le j<\lt(S)$, this condition is satisfied by at least one $\P\in\pset_S$.

\item Let $\alpha=Q_\j\in \bb_S$, as in (\ref{generic}). Then,
$$
w(\alpha)=H'_{1,j_1}+\cdots +H'_{r-1,j_{r-1}}+H_{r,j}-\lfloor H'_{1,j_1}+\cdots +H'_{r-1,j_{r-1}}+H_{r,j}\rfloor.
$$
In particular, $0\le w(\alpha)<1$. Moreover, for all $\P\in\pset_S$, we have
$w(\alpha)=w_\P(\alpha)$ if and only if $\psi_\P\nmid R_r(Q_{r,j})$.
\item Suppose that $\q\in\pset$ either belongs to $\pset_0$, or it does not dominate the prime ideals in $\pset_S$. Then, $w_\q(\alpha)>w(\alpha)$, for all $\alpha\in\bb_S$.
\end{enumerate}
\end{lemma}

\begin{proof}
Write $\lambda_r=-h_r/e_r$, with $h_r,e_r$ positive coprime integers.
Let $\P\in\pset_S$, and consider the type $\ty'_\P:=\ty_{\lambda_r,\psi_\P}=(\ty;(\phi_r,\lambda_r,\psi_\P))$, dividing $F_\P$. 

The shape of $N_r^-(Q_{r,j})$ is shown in Figure \ref{figNewQ}. The ordinate $H$ of the intersection point of the vertical axis with the line of slope $\lambda_r$ that  first touches $N_r^-(Q_{r,j})$ from below is equal to $y_{r,j}-(b_r-j)V_r$. 
By Proposition \ref{vgt} applied to the type $\ty'_\P$, 
$$
w_\P(Q_{r,j}(\t))\ge H/e_0\cdots e_{r-1}=H_{r,j},
$$ 
and equality holds if and only if $\psi_\P\nmid R_r(Q_{r,j})$.

Let $\varphi:=\prod_{\P\in\pset_S}\psi_\P$, so that $\lt(S)=e_r\deg\varphi$.
If $0\le j<\lt(S)$, Corollary \ref{RQ} shows that, 
$$
\deg R_r(Q_{r,j})\le j/e_r<\lt(S)/e_r=\deg\varphi.
$$
Since $\varphi$ is a separable polynomial, at least one irreducible factor $\psi_\P$, of $\varphi$, does not divide $R_r(Q_{r,j})$. This proves item 1.

Consider $\alpha=Q_\j\in\bb_S$, as in (\ref{generic}). Take arbitrary prime ideals $\q\in\pset$, $\P\in\pset_S$.
By the properties of the Okutsu frame (\ref{frame}), $w_\q(\t^{j_0})\ge 0=w_\P(\t^{j_0})$. By Theorem \ref{denquot} and Corollary \ref{applications}:
\begin{equation}\label{thm3.3}
w_\q(Q'_{i,j_i}(\t))\ge H'_{i,j_i}=w_\P(Q'_{i,j_i}(\t)),\quad \forall\,1\le i< r.
\end{equation}
By Theorem \ref{denquot}, $w_\q(Q_{r,j}(\t))\ge H_{r,j}$, and this coincides with $w_\P(Q_{r,j}(\t))$ if and only if 
 $\psi_\P\nmid R_r(Q_{r,j})$, by item 1. This proves item 2.

In order to prove item 3, it suffices to show that $w_\q(Q_{r,j}(\t))>H_{r,j}$, if $\q\in\pset_0$, or $\q\not\succ\P$.
Let us apply Theorem \ref{denquot} to the type $\ty$. If $\q$ falls in cases (ii) or (iii) of the proof of the theorem, the inequalities (\ref{stronger2}) and (\ref{stronger3}) show that $w_\q(Q_{r,j}(\t))>H_{r,j}$.  
Suppose that $\q$ falls in case (i); that is, $\ty\mid F_\q$. By (\ref{thmpol}) and Corollary \ref{previous},
\begin{equation}\label{values}
w_\q(\phi_r(\t))=(V_r+|\mu|)/e_0\cdots e_{r-1},\quad 
w_\P(\phi_r(\t))=(V_r+|\lambda_r|)/e_0\cdots e_{r-1},
\end{equation}
where $\mu$ is one of the slopes of $N_r^-(f)$. In the notation of Definition \ref{domination}, we have $\phi_S:=\phi_r$. Thus, $\q\not\succ\P$ means, by definition, $w_\q(\phi_r(\t))<w_\P(\phi_r(\t))$. By (\ref{values}), we get $|\mu|<|\lambda_r|$, and by (\ref{stronger1}), we deduce that $w_\q(Q_{r,j}(\t))>H_{r,j}$.
\end{proof}

\begin{definition}
We define a global $\star$-product on $B_\p\setminus\{0\}$ by:
$$
\alpha\star\beta:=\alpha\beta/\pi^{\lfloor w(\alpha\beta)\rfloor}\in B_\p,\quad \forall\,\alpha,\beta\in B_\p.
$$ 
It is clearly associative and commutative.
\end{definition}

\begin{lemma}\label{bsreduced}
Let $S$ be a terminal side. Let $V:=\prod_{\P\in\pset}\F_\P$, $V_S:=\prod_{\P\in\pset_S}\F_\P$, and $\op{pr}_S\colon V\lra V_S$ the canonical projection. Let $\bb_{S,\delta}:=\{\alpha\in\bb_S\mid w(\alpha)=\delta\}$, for some $\delta\in w(\bb_S)$. Then, $\op{pr}_S(\red_\delta(\bb_{S,\delta}))$ is an $\F_\p$-basis of $V_S$.
\end{lemma}
 
\begin{proof}
We keep the notation from Lemma \ref{auxiliar}. Let $\P\in\pset_S$, and denote $f_\P:=\deg\psi_\P$. The type $\ty'_\P=(\ty;(\phi_r,\lambda_r,\psi_\P))$ is $F_\P$-complete; hence \cite[Cor. 3.8]{HN},
\begin{equation}\label{ef}
e(\P/\p)=e_0\cdots e_r,\quad f(\P/\p)=f_0f_1\cdots f_{r-1}f_{\P}.
\end{equation}

For all $\P\in\pset_S$, the types $\ty'_\P$ coincide, except for the data involving
the polynomials $\psi_{\P}$. In particular, the tower of fields, $\F_\p=\F_0\subset\cdots\subset\F_r$, and the rational fractions $\pi_0,\dots,\pi_{r+1}\in K(x)$ of (\ref{ratfracs}), are the same for all $\P\in\pset_S$. We denote:
$$
\F_\P=\F_r[y]/(\psi_{\P}(y))=\F_r[z_{\P}],\quad \Pi:=\Pi_\P=\pi_{r+1}(\t),\quad \Pi_0:=\pi_{r}(\t),
$$
where $z_\P$ is the class of $y$ in $\F_\P$, so that $\psi_{\P}(z_\P)=0$. Recall that $w_\P(\Pi)=1/(e_0\cdots e_r)$, $w_\P(\Pi_0)=1/(e_0\cdots e_{r-1})$, for all $\P\in\pset_S$.

Consider the set:
$$
\bb_S^0:=\left\{\dfrac{\t^{j_0} Q'_{1,j_1}(\t)\cdots Q'_{r-1,j_{r-1}}(\t)}{\pi^{\lfloor H'_{1,j_1}+\cdots H'_{r-1,j_{r-1}}\rfloor}}\ \Big|\ 0\le j_i<e_if_i,\ 0\le i<r\right\}.
$$
For all $\P\in\pset_S$, we have:
\begin{enumerate}
\item[(i)] $i_\P(\bb_S^0)$ is a level $r-1$ basis in standard form of $L_\P/K_\p$,
\item[(ii)] $w(\bb_S^0)=w_\P(\bb_S^0)=\{t'/(e_0\cdots e_{r-1})\mid t'\in\Z,\ 0\le t'<e_0\cdots e_{r-1}\}$. 
\end{enumerate}
In fact, we saw (i) along the proof of Theorem \ref{Pbasis2}, and (ii) is deduced from (\ref{thm3.3}). 

Let $\delta\in w(\bb_S)$. By Lemma \ref{auxiliar}, $0\le\delta<1$ and there exists $\P\in\pset_S$ such that $\delta\in w_{\P}(L)=(e_0\cdots e_r)^{-1}\Z$. Thus, $\delta=t/(e_0\cdots e_r)$, for some integer $0\le t< e_0\cdots e_r$. Let $q_t$ be the unique integer, $0\le q_t<e_r$, such that $q_th_r\equiv t\md{e_r}$. 
For any $0\le j<\lt(S)$, the argument of the proof of item (a) of Lemma \ref{valuett} shows that:
$$
j\not\equiv q_t\md{e_r}\ \Longrightarrow\ H_{r,j}+t'/(e_0\cdots e_{r-1})\not\equiv\delta\md{\Z},
$$
for all integers, $0\le t'<e_0\cdots e_{r-1}$, whereas
$$
j=q_t+k e_r\ \Longrightarrow\ H_{r,j}+t_k/(e_0\cdots e_{r-1})\equiv\delta\md{\Z},
$$
for a uniquely determined integer $0\le t_k<e_0\cdots e_{r-1}$. 
This leads to: 
$$
\bb_{S,\delta}=\bigcup\nolimits_{0\le k<f_S}Q_{r,q_t+ke_r}(\t)\star \bb_{S,t_k}^0,
$$
where $f_S:=\sum_{\P\in\pset_S}f_{\P}=\lt(S)/e_r=\dim_{\F_r}V_S$, and  $\bb^0_{S,t_k}$ is the subset of $\bb_S^0$ formed by those $\alpha^0$ such that $w(\alpha^0)=t_k/e_0\cdots e_{r-1}$. 

For any $0\le k<f_S$, write $\bb_{S,t_k}^0=\Pi_0^{t_k}U_k$, for $U_k\subset L$. By condition (i) above, $\red_{L_\P}(i_\P(U_k))\subset \F_r$ is an $\F_\p$-basis of $\F_r$, for all $\P\in\pset_S$. 
Now, the elements in $\bb_{S,\delta}$ may be parameterized as: 
$$\alpha_{k,u}=Q_{r,j}(\t)\star\Pi_0^{t_k}u,\quad 0\le k<f_S,\ u\in U_k,
$$for $j=q_t+ke_r$. Let us compute the $\P$-th component, $\red_{L_\P}(i_\P(\alpha_{k,u}/\Pi^t))$, of $\red_\delta(\alpha_{k,u})\in V$, for all $\P\in\pset_S$. 
By item 2 of Lemma \ref{auxiliar}:
\begin{equation}\label{recall}
w(\alpha_{k,u})=w_\P(\alpha_{k,u})\sii \psi_\P\nmid R_r(Q_{r,j})\sii R_r(Q_{r,j})(z_\P)\ne0. 
\end{equation} 

If $\P$ satisfies (\ref{recall}), then $i_\P$ is compatible with the global and local $\star$ operations, and the arguments of the proof of Proposition \ref{PbasisFinal}, lead to (\ref{epsk2}), (\ref{epsk3}), and:
\begin{align*}
\eta_{k,u,\P}&:= \red_{L_\P}(i_\P(\alpha_{k,u}/\Pi^t))=\red_{L_\P}(i_\P(Q_{r,j}(\t)\star\Pi_0^{t_k}u\,\Pi^{-t}))\\&=\red_{L_\P}(i_\P(Q_{r,j}(\t))\star i_\P(\Pi_0)^{t_k}i_\P(\Pi)^{-t})\,\bar{u}\\&=R_r(Q_{r,j})(z_\P)\cdot\tau_k\cdot (z_\P)^{N+k-d_k}\,\bar{u}\\&=\zeta_k\cdot (z_\P)^N(C_{d_k}(z_\P)^{k-d_k}+\cdots+C_0(z_\P)^k)\,\bar{u},
\end{align*}
where $\bar{u}:=\red_{L_\P}(i_\P(u))$, $d_k=\deg R_r(Q_{r,j})\le k$, $\tau_k,\zeta_k,C_0,C_{d_k}\in\F_r^*$, $C_i\in\F_r$, for $i\ne0,d_k$, and $N$ is an integer that depends only on $t$ (that is, on $\delta$). 

If $\P$ does not satisfy (\ref{recall}), then $w_\P(\alpha_{k,u})>w(\alpha_{k,u})$, so that $\eta_{k,u,\P}=0$. Since $R_r(Q_{r,j})(z_\P)=0$ in this case, the above formula for $\eta_{k,u,\P}$ holds for all $\P\in\pset_S$.   

Our aim is to show that the vectors $\eta_{k,u}:=(\eta_{k,u,\P})_{\P\in\pset_S}\in V_S$, for $0\le k<f_S$ and $u\in U_k$, are an $\F_\p$-basis of $V_S$. Clearly, the map:
$$
(x_\P)_{\P\in\pset_S}\mapsto ((C_0)^{-1}(z_\P)^{-N}x_\P)_{\P\in\pset_S}
$$
is an $\F_\p$-automorphism of $V_S$. Thus, since $\bar{u},\zeta_k\in \F_r^*$ do not depend on $\P\in\pset_S$, we may assume that
$$
\eta_{k,u}=\bar u\cdot\zeta_k\cdot \eta_k,\quad \eta_k=\left(c_{d_k}(z_\P)^{k-d_k}+\cdots+c_1(z_\P)^{k-1}+(z_\P)^k\right)_{\P\in\pset_S}\in V_S,
$$
where $c_i:=C_i/C_0$, for all $i$. Since for all $k$ the family $\{\bar{u}\mid u\in U_k\}$ is an $\F_\p$-basis of $\F_r$, it suffices to check that the family of all $\{\zeta_k\eta_k\mid 0\le k<f_S\}$ is an $\F_r$-basis of $V_S$. Since $\dim_{\F_r}V_S=f_S$ and all $\zeta_k$ belong to $\F_r^*$, this is equivalent to $\{\eta_k\mid 0\le k<f_S\}$ being an $\F_r$-linearly independent family of $V_S$. We can relate this family to the family $\eta'_k:=((z_\P)^k)_{\P\in\pset_S}$ by the following equations:
$$
\eta_k=\eta'_k+c_1\eta'_{k-1}+\cdots+c_{d_k}\eta'_{k-d_k}.
$$ 
Since the transition matrix between the two families is invertible, it suffices to check that the family $\{\eta'_k\mid 0\le k<f_S\}$ is $\F_r$-linearly independent. Now, a linear relation of the form: $\sum_{0\le k<f_S}a_k\eta'_k=0$, with $a_k\in \F_r$, is equivalent to: 
$$
\sum\nolimits_{0\le k<f_S}a_k(z_\P)^k=0,\quad\forall\,\P\in\pset_S.
$$  
Since the irreducible polynomials $\psi_\P$, for $\P\in\pset_S$, are pairwise different, this implies that the polynomial $\sum_{0\le k<f_S}a_kx^k$ is divisible by the polynomial $\prod_{\P\in\pset_S}\psi_\P$, which has degree $f_S$. This occurs only when all coefficients $a_k$ vanish. 
\end{proof}

\noindent{\sl Proof of Theorem \ref{pBasis2}.} 
Denote $n_\P=e(\P/\p)f(\P/\p)$, for all $\P\in\pset$. Clearly, 
$$
\#\bb_{\pset_0}=\sum\nolimits_{\P\in\pset_0}f_{0,\P}=\sum\nolimits_{\P\in\pset_0}n_\P.
$$
On the other hand, by (\ref{ef}), for any $S\in\st$, 
\begin{align*}
\#\bb_S&=(e_0f_0)\cdots (e_{r-1}f_{r-1})\,\lt(S)\\
&=(e_0f_0)\cdots (e_{r-1}f_{r-1})\cdot e_r\cdot \sum\nolimits_{\P\in\pset_S}f_\P=\sum\nolimits_{\P\in\pset_S}n_\P.
\end{align*}
Thus, by Lemma \ref{partition}, $\#\bb=\sum_{\P\in\pset}n_\P=n$.
Also, $\bb\subset B_\p$, by construction. 

For any $\P\in\pset_0$, Corollary \ref{applications} shows that $w_\P(Q_\P(\t))=0$, and (\ref{frame}) shows that $w_\P(\t^j)=0$, for all $0\le j<f_{0,\P}$; hence, $w(\alpha)=0$ for all $\alpha\in\bb_{\pset_0}$. We conclude that $0\le w(\alpha)<1$, for all $\alpha\in\bb$, by
item 2 of Lemma \ref{auxiliar}. Therefore, by Lemma \ref{reducedbasis}, we need only to check that the set $\bb$ is reduced. To this end, we apply the criterion of Lemma \ref{criterion}.

For any $\delta\in w(\bb)$, let $\bb_{S,\delta}:=\bb_\delta\cap \bb_S$. The set $\bb_\delta$ splits into the disjoint union:
$$
\bb_\delta=\left\{
\begin{array}{ll}
 \bb_{\pset_0}\cup\left(\bigcup_{S\in\st}\bb_{S,\delta}\right),&\mbox{ if }\delta=0,\\ 
\bigcup_{S\in\st}\bb_{S,\delta},&\mbox{ if }\delta>0. 
\end{array}\right.$$

Our aim is to prove the $\F_\p$-linear independence of the family
$\red_\delta(\bb_\delta)$, for all $\delta\in w(\bb)$. Let us show first that the family $\bigcup_{S\in\st}\red_\delta\left(\bb_{S,\delta}\right)$ is linearly independent. 

Take any $\delta\in w(\bb)$. Let $\mathcal{S}_\delta:=\{S\in\st\mid \bb_{S,\delta}\ne\emptyset\}$. For any $S\in\mathcal{S}_\delta$, write: $\red_\delta(\bb_{S,\delta})=\{\zeta_{m,S}\mid 1\le m\le n_{S,\delta}\}\subset V$. 

Suppose that for some family of elements $a_{m,S}\in\F_\p$, we have
\begin{equation}\label{equalzero}
\sum\nolimits_{m,S}a_{m,S}\,\zeta_{m,S}=0,
\end{equation}
the sum running on $S\in\mathcal{S}_\delta$ and $1\le m\le n_{S,\delta}$. Take $T\in\mathcal{S}_\delta$ minimal with respect to the relationship of domination (cf. the remark following Lemma \ref{partorder}); that is:
$$
\q\not\succ\P,\ \forall\,\q\in\pset_{T},\ \forall\,\P\in\pset_S,\ \forall\,S\in \mathcal{S}_\delta,\ S\ne T.
$$
By item 3 of Lemma \ref{auxiliar}, 
$$
w_\q(\alpha)>\delta,\quad \forall\,\q\in\pset_{T},\ \forall\,\alpha\in\bb_{S,\delta},\ \forall\,S\in\mathcal{S}_\delta, \ S\ne T. 
$$
Hence, $\op{pr}_{T}(\zeta_{m,S})=0$, for all $S\in\mathcal{S}_\delta$, $S\ne T$, and all $m$. Thus, if we apply $\op{pr}_{T}$ to both sides of (\ref{equalzero}), we get
$$
\sum\nolimits_ma_{m,T}\op{pr}_{T}(\zeta_{m,T})=0,
$$
and by Lemma \ref{bsreduced}, $a_{m,T}=0$, for all $m$. Thus, we get again an equation like (\ref{equalzero}), for $S$ running on the set $\mathcal{S}_\delta\setminus\{T\}$. By applying in a recurrent way the same argument, we conclude that $a_{m,S}=0$ for all $m,S$. Therefore, the family $\bigcup_{S\in\st}\red_\delta\left(\bb_{S,\delta}\right)$ is $\F_\p$-linearly independent.

This proves the theorem in the case $\delta>0$. Suppose now $\delta=0$.

Define the \emph{support} of a vector $(x_\P)_{\P\in\pset}\in V$, as the set of indices $\P\in\pset$ such that $x_\P\ne 0$. For each $\P\in\pset_0$,  $\red_0(\bb_\P)$ is an $\F_\p$-linearly independent subset of $V$, and all these vectors have support $\{\P\}$, because $w_\q(Q_\P)>0$, for all $\q\in\pset$, $\q\ne\P$. In fact, since $\psi_{0,\q}\ne\psi_{0,\P}$, the prime ideal $\q$ falls in case (iii) of Theorem \ref{denquot}, applied to the type of order zero $\ty=(\psi_{0,\P})$, with representative $\phi_{1,\P}$. 

We have seen that the subset $\bigcup_{S\in\st}\red_0\left(\bb_{S,0}\right)$ is $\F_\p$-linearly independent.  All these vectors in $V$ have support contained in $\pset\setminus\pset_0$, by item 3 of Lemma \ref{auxiliar}. Therefore, $\red_0(\bb_0)$ is $\F_\p$-linearly independent, because it is the union of linearly independent subsets with pairwise disjoint supports.\qed

\subsection{Complexity analysis} 
Denote $\delta:=v_\p(\dsc(f))$. If $\delta=0$, then $B_\p=A_\p[\t]$; thus, we assume $\delta>0$ in our analysis.
By \cite[Thm. 3.14]{BNS}, for the computation of a $\p$-integral basis of $B/A$ we may work modulo $\p^{\delta+1}$. Hence, we assume that the elements of $A$ are finite $\pi$-adic developments of length $\delta+1$. 

\begin{definition}
An operation in $A$ is called \emph{$\p$-small} if it involves two elements belonging to a fixed system of representatives of $A/\p$.
\end{definition}

Each multiplication in $A$ costs $O(\delta^{1+\epsilon})$ $\p$-small operations, if we assume the fast multiplications techniques of Sch\"onhage-Strassen \cite{SS}. Also, if $q:=\#A/\p$, a $\p$-small operation in $A$ requires  $O\left(\log(q)^{1+\epsilon}\right)$ word operations, the cost of an operation in the residue field $A/\p$. 

The Montes algorithm has a cost of $O\left(n^{2+\epsilon}+n^{1+\epsilon}\delta\log q+n^{1+\epsilon}\delta^{2+\epsilon}\right)$
$\p$-small operations \cite[Thm. 5.15]{BNS}. Let us estimate the cost of the extra tasks that are necessary to compute the $\p$-integral basis.
The computation of $\bb_{\pset_0}$ being negligible, let us discuss the computation of $\bigcup_{S\in\st}\bb_S$.

For any element $\alpha\in\bb_S$, the factors $Q_{i,j}$ of the numerator of $\alpha$, and the exponents $H_{i,j}$ of the denominator, are computed along the flow of the Montes algorithm. The final computation of $\alpha$ is dominated by the product of the numerators in $A[\t]$, and we may neglect the computation of the sum of the $H_{i,j}$ and the final division by a power of $\pi$.  

Suppose $S$ is a terminal side of order $r$, derived from a type $\ty$ of order $r-1$, with representative $\phi_r$. We  compute $\bb_S$ in $r$ steps:
$$
\bb_{S,0}=\{1,\t,\dots,\t^{f_0-1}\},\quad \bb_{S,i}=\bigcup\nolimits_{0\le j<e_if_i}Q'_{i,j}(\t)\star\bb_{S,i-1},\ 1\le i<r,
$$
and finally, $\bb_S=\bb_{S,r}=\bigcup_{0\le j<e_rf_S}Q_{r,j}(\t)\star\bb_{S,r-1}$. Clearly,
$$
\begin{array}{l}
 \#\bb_{S,i}=(e_0f_0)\cdots (e_if_i),\quad 0\le i<r,\\
n_S:=\#\bb_{S,r}=(e_0f_0)\cdots (e_{r-1}f_{r-1})e_rf_S=\sum_{\P\in\pset_S}n_\P. 
\end{array}
$$

If we keep only the numerators in mind, each element of $\bb_{S,i}$, $i\ge 1$, is obtained after one multiplication in $A[\t]$: $Q'_{i,j}(\t)$, or $Q_{r,j}(\t)$, times an element in $\bb_{S,i-1}$. Thus, the total number of multiplications for the computation of $\bb_S$ is:
$$
N=e_0f_0+(e_0f_0)(e_1f_1)+\cdots +(e_0f_0)\cdots (e_{r-1}f_{r-1})+(e_0f_0)\cdots (e_{r-1}f_{r-1})(e_rf_S).
$$
Now, since the type $\ty$ is optimal, we have $e_if_i\ge 2$, for $1\le i<r$. Thus,
$$
N\le \dfrac{n_S}{2^{r-1}}+\dfrac{n_S}{2^{r-2}}+\cdots+ \dfrac{n_S}{1}+n_S=n_S(1+2-2^{1-r})\le 3n_S.
$$
Since $n\ge\sum_{S\in\st}n_S$, the total number of multiplications in $A[\t]$ required for the computation of $\bb$ is $O(n)$. Each multiplication in $A[\t]$ has a cost of $O(n^{1+\epsilon})$ multiplications in $A$. Therefore, the total cost of the computation of $\bb$ is $O\left(n^{2+\epsilon}\delta^{1+\epsilon}\right)$ $\p$-small operations. Adding to this cost the cost of the Montes algorithm, we get the following total estimation.
 
\begin{thm}\label{complexity}
The computation of a $\p$-integral basis of $B/A$ requires not more than $O\left(n^{2+\epsilon}\delta^{1+\epsilon}+n^{1+\epsilon}\delta\log q+n^{1+\epsilon}\delta^{2+\epsilon}\right)$ $\p$-small operations in $A$.  If $A/\p$ is small, we obtain an estimation of $O\left(n^{2+\epsilon}\delta^{1+\epsilon}+n^{1+\epsilon}\delta^{2+\epsilon}\right)$ word operations. 
\end{thm}

\subsection{An example}
Let us show how the method of the quotients works with an example. Take $A=\Z$, $\p=2\Z$ and $$f(x)=x^{12} + 14x^{10} + 60x^8 + 32x^7 + 80x^6 + 128x^5 - 80x^4 + 256x^3 - 288x^2 - 256x +832.$$
Since $f(x)\equiv x^{12}\md2$, there will be only one tree of types, whose root node is the type of order zero, $\ty_0=(y)$, determined by the irreducible polynomial $\psi_0(y)=y$. Take $\phi_1(x)=x$ as a representative of this type. The Newton polygon $N_1(f)$ is one-sided, and it has slope $-1/2$.

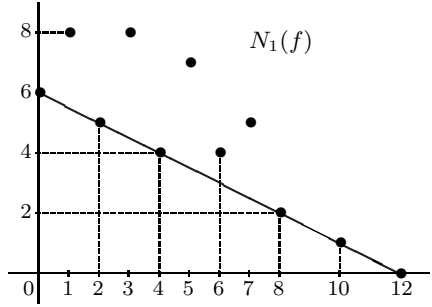
\begin{center}
\setlength{\unitlength}{4.mm}
\begin{picture}(14,10.5)
\put(-.15,6.8){$\bullet$}\put(.85,8.8){$\bullet$}
\put(1.85,5.8){$\bullet$}\put(2.85,8.8){$\bullet$}
\put(3.85,4.8){$\bullet$}\put(4.85,7.8){$\bullet$}
\put(5.85,4.8){$\bullet$}\put(6.85,5.8){$\bullet$}
\put(7.85,2.8){$\bullet$}\put(9.85,1.8){$\bullet$}
\put(11.85,.8){$\bullet$}
\put(0,0){\line(0,1){10}}\put(-1,1){\line(1,0){14}}
\put(0,7.03){\line(2,-1){12}}\put(0,7){\line(2,-1){12}}
\put(7,8.5){\begin{small}$N_1(f)$\end{small}}
\put(11.7,.3){\begin{footnotesize}$12$\end{footnotesize}}
\put(9.6,.3){\begin{footnotesize}$10$\end{footnotesize}}
\put(7.8,.3){\begin{footnotesize}$8$\end{footnotesize}}
\put(6.8,.3){\begin{footnotesize}$7$\end{footnotesize}}
\put(5.8,.3){\begin{footnotesize}$6$\end{footnotesize}}
\put(4.8,.3){\begin{footnotesize}$5$\end{footnotesize}}
\put(3.8,.3){\begin{footnotesize}$4$\end{footnotesize}}
\put(2.8,.3){\begin{footnotesize}$3$\end{footnotesize}}
\put(1.8,.3){\begin{footnotesize}$2$\end{footnotesize}}
\put(.8,.3){\begin{footnotesize}$1$\end{footnotesize}}
\put(-.5,.3){\begin{footnotesize}$0$\end{footnotesize}}
\put(-.6,8.85){\begin{footnotesize}$8$\end{footnotesize}}
\put(-.6,6.85){\begin{footnotesize}$6$\end{footnotesize}}
\put(-.6,4.85){\begin{footnotesize}$4$\end{footnotesize}}
\put(-.6,2.85){\begin{footnotesize}$2$\end{footnotesize}}
\put(1,.9){\line(0,1){.2}}\put(2,.9){\line(0,1){.2}}
\put(3,.9){\line(0,1){.2}}\put(4,.9){\line(0,1){.2}}
\put(5,.9){\line(0,1){.2}}\put(6,.9){\line(0,1){.2}}
\put(7,.9){\line(0,1){.2}}\put(8,.9){\line(0,1){.2}}
\put(10,.9){\line(0,1){.2}}
\multiput(2,.9)(0,.25){21}{\vrule height2pt}
\multiput(4,.9)(0,.25){17}{\vrule height2pt}
\multiput(6,.9)(0,.25){17}{\vrule height2pt}
\multiput(8,.9)(0,.25){9}{\vrule height2pt}
\multiput(10,.9)(0,.25){5}{\vrule height2pt}
\multiput(-.1,9)(.25,0){5}{\hbox to 2pt{\hrulefill }}
\multiput(-.1,5)(.25,0){17}{\hbox to 2pt{\hrulefill }}
\multiput(-.1,3)(.25,0){33}{\hbox to 2pt{\hrulefill }}
\end{picture}
\end{center}\bigskip

 The residual polynomial of the first order attached to the side $S=N_1(f)$ is:
$$
R_{-1/2,1}(f)(y)=y^6+y^5+y^4+y^2+y+1=(y^2+y+1)(y+1)^4.
$$
Thus, the type $\ty_0$ ramifies into two types of order one:
$$
\ty_1=(y;(x,-1/2,y^2+y+1)),\quad \ty'_1=(y;(x,-1/2,y+1)).
$$ 
Since $y^2+y+1$ divides $R_{-1/2,1}(f)(y)$ with exponent one, the type $\ty_1$ is $f$-complete and $S$ is a terminal side of order $1$, with $\lt(S)=4$. The type $\ty_1$ singles out a prime ideal $\P$ with $e(\P/p)=f(\P/p)=2$. 
If $Q_1,\dots,Q_{12}$ are the twelve quotients of the $x$-adic development of $f(x)$, we have: 
$$
\begin{array}{ll}
Q_{1,0}=Q_{12}=1, \quad H_{1,0}=0,&\quad Q_{1,1}=Q_{11}=x,\quad H_{1,1}=1/2,\\
Q_{1,2}=Q_{10}=x^2+14, \quad H_{1,2}=1,&\quad Q_{1,3}=Q_9=x^3+14x,\quad H_{1,3}=3/2.
\end{array}
$$
The set $\bb_S$ contains the following four globally integral elements:
$$
\bb_S=\left\{1,\ \t,\ (\t^2+14)/2,\ (\t^3+14\t)/2\right\}.
$$

The type $\ty'_1$ is not complete, and its analysis requires some more work in order two. Before analyzing its branching, we store a list $\bb_1=\{1,\t\}$, with the quotients $Q'_{1,j}$, for $0\le j<e_1f_1=2$, and also the corresponding values $H'_{1,0}=0$, $H'_{1,1}=1/2$. All future branches of $\ty'_1$ will share these data.
 
Let us choose $\phi_2(x)=x^2+2$ as a representative of $\ty_1'$. By Lemma \ref{length}, $\ell(N_{2}^-(f))=\ord_{\ty_1'}(f)=4$; thus, we compute the $\phi_2$-expansion of $f(x)$ only up to degree four:
$$
f(x)= 1024-512x+ 128x\phi_2(x)-64x\phi_2(x)^2+32x\phi_2(x)^3-20\phi_2(x)^4+....
$$
The first four quotients of this $\phi_2$-development are:
$$
\begin{array}{l}
Q_1=x^{10} + 12x^8 + 36x^6 + 32x^5 + 8x^4 + 64x^3 - 96x^2 + 128x - 96,\\
Q_2=x^8 + 10x^6 + 16x^4 + 32x^3 - 24x^2 - 48,\\
Q_3=x^6 + 8x^4 + 32x - 24,\\
Q_4=x^4 + 6x^2 - 12.
\end{array}
$$
Since $V_2=v_2(\phi_2)=v_2(2)=2$, and $v_2(x)=1$, we get 
$v_2(1024-512x)=19$, $v_2(128x\phi_2)=v_2(64x(\phi_2)^2)=v_2(32x(\phi_2)^3)=17$, $v_2(20(\phi_2)^4)=12$. Therefore, $N_2^-(f)$ has two sides of slopes $-2$ and $-5/3$:\bigskip

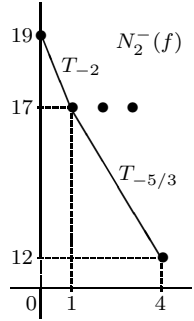
\begin{center}
\setlength{\unitlength}{4.mm}
\begin{picture}(8,10)
\put(-.2,9.2){$\bullet$}\put(.85,6.8){$\bullet$}
\put(1.85,6.8){$\bullet$}\put(2.85,6.8){$\bullet$}
\put(3.85,1.8){$\bullet$}
\put(0,0){\line(0,1){1}}\put(0,2){\line(0,1){8.5}}
\put(-1,1){\line(1,0){6}}
\put(1,7.03){\line(-2,5){1}}\put(1,7){\line(-2,5){1}}
\put(1,7){\line(3,-5){3}}\put(1,7.03){\line(3,-5){3}}
\put(.7,8.2){\begin{footnotesize}$T_{-2}$\end{footnotesize}}
\put(2.6,4.6){\begin{footnotesize}$T_{-5/3}$\end{footnotesize}}
\put(2.5,9){\begin{footnotesize}$N_2^-(f)$\end{footnotesize}}
\multiput(0,1)(0,.25){4}{\vrule height2pt}
\multiput(4,.9)(0,.25){5}{\vrule height2pt}
\multiput(1,.9)(0,.25){25}{\vrule height2pt}
\put(3.8,.3){\begin{footnotesize}$4$\end{footnotesize}}
\put(.85,.3){\begin{footnotesize}$1$\end{footnotesize}}
\put(-.5,.3){\begin{footnotesize}$0$\end{footnotesize}}
\put(-1,9.2){\begin{footnotesize}$19$\end{footnotesize}}
\put(-1,6.8){\begin{footnotesize}$17$\end{footnotesize}}
\put(-1,1.8){\begin{footnotesize}$12$\end{footnotesize}}
\multiput(-.1,2)(.25,0){17}{\hbox to 2pt{\hrulefill }}
\multiput(-.1,7)(.25,0){5}{\hbox to 2pt{\hrulefill }}
\put(1,.9){\line(0,1){.2}}\put(4,.9){\line(0,1){.2}}
\end{picture}
\end{center}\medskip

Since both sides have degree one, the residual polynomials of second order have degree one: $R_{-2,2}(f)(y)=R_{-5/3,2}(f)(y)=y+1$. Thus, the type $\ty_1'$ branches into two $f$-complete types of order two:
$$
\ty_2=(y;(x,-1/2,y+1);(\phi_2,-2,y+1)),\quad \ty_2'=(y;(x,-1/2,y+1);(\phi_2,-5/3,y+1)).
$$ They single out two prime ideals $\q$, $\q'$,  with $e(\q/2)=2$, $e(\q'/2)=6$, $f(\q/2)=f(\q'/2)=1$. 

Both sides $T_{-2}$ and $T_{-5/3}$ are terminal. With respect to $T_{-2}$, we have $Q_{2,0}=Q_1$, $H_{2,0}=15/2$, so that
$$
\bb_{T_{-2}}=Q_{2,0}\star \bb_1=\left\{Q_1(\t)/2^7,\,\t Q_1(\t)/2^8\right\}.
$$
With respect to $T_{-5/3}$, we have $Q_{2,0}=Q_4$, $H_{2,0}=2$, $Q_{2,1}=Q_3$, $H_{2,1}=23/6$, $Q_{2,2}=Q_2$, $H_{2,2}=17/3$, so that
$$
\bb_{T_{-5/3}}=\bigcup_{0\le j<3}Q_{2,j}\star \bb_1=\left\{\dfrac{Q_4(\t)}{2^2},\,\dfrac{\t Q_4(\t)}{2^2},\,\dfrac{Q_3(\t)}{2^3},\,\dfrac{\t Q_3(\t)}{2^4},\,\dfrac{Q_2(\t)}{2^5},\,\dfrac{\t Q_2(\t)}{2^6}\right\}.
$$
The $2$-integral basis, $\bb=\bb_S\cup\bb_{T_{-2}}\cup\bb_{T_{-5/3}}$, is complete.
The Hermite Normal Form algorithm transforms the basis into the following $12$ integral elements:
$$
\as{2.4}
\begin{array}{l}
1, \,\t,\,\dfrac{\t^2}2,\, \dfrac{\t^3}{2},\, \dfrac{\t^4}{4},\, \dfrac{\t^5}{4}, \dfrac{\t^6}{8},\,\dfrac{\t^7+8\t}{16},\
\dfrac{\t^8+2\t^6+8\t^2+16}{32},\,\dfrac{\t^9+2\t^7+8\t^3+16\t}{64},\\,\dfrac{\t^{10}+12\t^6+8\t^4+96}{128},\, \dfrac{\t^{11}+12\t^7+8\t^5+64\t^3+224\t}{256},
\end{array}
$$

\end{document}